\documentclass[10pt]{amsart}

\newtheorem{proposition}{Proposition}[section]
\newtheorem{lemma}[proposition]{Lemma}

\newtheorem{theorem}[proposition]{Theorem}

\theoremstyle{definition}

\theoremstyle{remark}
\newtheorem{remark}[proposition]{Remark}

\newcommand{\selabel}[1]{\label{se:#1}}
\newcommand{\seref}[1]{Section~\ref{se:#1}}

\newcommand{\eqlabel}[1]{\label{eq:#1}}

\newcommand{\Hom}{{\rm Hom}}

\newcommand{\ev}{{\rm ev}}
\newcommand{\coev}{{\rm coev}}

\def\lan{\langle}
\def\ran{\rangle}
\def\ot{\otimes}

\def\cal{\mathcal}

\newcommand{\Cc}{\mathcal{C}}

\newcommand{\Mm}{\mathcal{M}}

\def\*C{{}^*\hspace*{-1pt}{\Cc}}
\def\text#1{{\rm {\rm #1}}}

\def\yd{\mbox{$_H^H{\mathcal YD}$}}
\def\va{\varepsilon}
\def\v{\varphi}

\def\tr{\triangleright}
\def\rh{\rightharpoonup}
\def\lh{\leftharpoonup}

\def\ra{\rightarrow}
\def\a{\alpha}
\def\b{\beta}

\def\l{\lambda}
\def\r{\rho}
\def\cd{\cdot}
\def\d{\delta}
\def\ov{\overline}
\def\un{\underline}
\def\mf{\mathfrak}

\newcommand{\smi}{\mbox{$S^{-1}$}}

\def\rawo\lonra{\longrightarrow}

\def\ot{\otimes}
\def\act{\hspace*{-1pt}\ra\hspace*{-1pt}}
\def\<{\langle}
\def\>{\rangle}
\newcommand{\tpla}{\mbox{$\tilde {p}^1$}}
\newcommand{\tplb}{\mbox{$\tilde {p}^2$}}

\newcommand{\tqla}{\mbox{$\tilde {q}^1$}}
\newcommand{\tqlb}{\mbox{$\tilde {q}^2$}}
\newcommand{\tQla}{\mbox{$\tilde {Q}^1$}}
\newcommand{\tQlb}{\mbox{$\tilde {Q}^2$}}

\begin{document}
\title[The representation-theoretic rank]{The representation-theoretic
rank of the doubles of quasi-quantum groups}
\author{Daniel Bulacu}
\address{Faculty of Mathematics and Informatics, University of Bucharest,
Str. Academiei 14, RO-010014, Bucharest 1, Romania}
\email{dbulacu@al.math.unibuc.ro}

\author{Blas Torrecillas}
\address{Department of Algebra and Analysis,
University of Almeria, 04071 Almeria, Spain}
\email{btorreci@ual.es}

\thanks{During the time of the preparation of this work the first author was
financially supported by the LIEGRITS program,
a Marie Curie Research Training Network funded by the European
community as project MRTN-CT 2003-505078.
He would like to thank to the guest, University of Almeria (Spain),
for their warm hospitality.}


\begin{abstract}
We compute the representation-theoretic
rank of a finite dimensional quasi-Hopf algebra $H$ and
of its quantum double $D(H)$, within the rigid braided category
of finite dimensional left $D(H)$-modules.
\end{abstract}
\maketitle
\section{Introduction}\selabel{1}
${\;\;\;}$ 
The definition of a quasi-bialgebra $H$ ensures that
the category of left $H$-modules ${}_H{\cal M}$ is a monoidal
category, and for a quasi-Hopf algebra $H$ the definition ensures
that ${}_H{\cal M}^{\rm fd}$, the category of finite dimensional
left $H$-modules, is a monoidal category with duality. Moreover, a
quasi-Hopf algebra is called quasi-triangular (ribbon) if the
monoidal category ${}_H{\cal M}$ is braided (ribbon, at least in
the finite dimensional case). So, in general, the study of
quasi-Hopf algebras is strictly connected to the study of
monoidal, or braided (ribbon) categories. Consequently, when we
want to define some classes of quasi-Hopf algebras the first thing
we should think about is to reword at a categorical level the
corresponding definitions given in the classical Hopf case. If it
is possible, then we can come back to the quasi-Hopf case. For
example, this was the case in \cite{bt}, where using the
categorical interpretation of a factorizable Hopf algebra (due to
Majid \cite{maj}), we were able to define and study the class of
factorizable quasi-Hopf algebras. But sometimes this point of view
cannot be followed. For further use we choose as an example the cosemisimple notion.\\
${\;\;\;}$ 
It is well known that a Hopf algebra $H$ is
cosemisimple if the category of left or right $H$-comodules is
cosemisimple. In the quasi-Hopf case we cannot consider
$H$-comodules, because the quasi-Hopf algebra $H$ is not
coassociative; thus in this case we have to look at some other
objects. One of these objects could be the quantum double $D(H)$
associated to a finite dimensional quasi-Hopf algebra $H$. In the
Hopf case we know that $D(H)$ is semisimple if and only if $H$ is
semisimple and cosemisimple (see \cite{r}). But, once again, at
this moment we cannot follow this path because in the quasi-Hopf
case we do not know the form of an integral in $D(H)$. However, in
the Hopf case, the Maschke-type theorem asserts that $H$ is
cosemisimple if and only if there exists a left or right integral
$\l$ in $H^*$ such that $\l (1)=1$. Now by \cite{ab} this is
equivalent to the existence of a bilinear form $\sigma \in (H\ot
H)^*$ such that $h_1\sigma (h_2, h')=\sigma (h, h'_1)h'_2$ and
$\sigma (h_1, h_2)=\va (h)$, for all $h, h'\in H$. This approach
was used by Hausser and Nill in \cite{hn3} for the quasi-Hopf
algebra setting. They proved that for a finite dimensional
quasi-Hopf algebra $H$ there is a one-to-one correspondence
between left cointegrals $\l \in H^*$ (see the definition below)
satisfying the normalized condition $\l (\smi (\a)\b)=1$ (here
$\a$ and $\b$ are the elements which occur in the definition of
the antipode $S$ of $H$), and certain bilinear forms $\sigma \in
(H\ot H)^*$ satisfying properties which generalize the ones
described above. This is why we will say that a finite dimensional
quasi-Hopf algebra $H$ is cosemisimple if $H$ admits a left
cointegral $\l$ obeying $\l (\smi (\a)\b)=1$. Furthermore, we believe
that an integral in $D(H)$ has the form $\b\rh \l \bowtie r$, so
if it is the case then $D(H)$ is semisimple if and only if $H$ is
semisimple and the left cointegral $\l$ satisfies $\l (\smi
(\a)\b)=1$ (here $r$ is a right integral in $H$). Comparing this
with the Hopf algebra case we will land to the same definition for
a finite dimensional
cosemisimple quasi-Hopf algebra.\\
${\;\;\;}$ 
The starting point of this paper was the intention to
generalize some important results concerning semisimple
cosemisimple Hopf algebras to quasi-Hopf algebras. Namely, a Hopf
algebra over a field of characteristic zero is semisimple if and
only if it is cosemisimple, if and only if it is involutory, this
means $S^2=id_H$. The result was proved by Larson and Radford in
\cite{lr1, lr2}, answering in positive, in characteristic zero,
the fifth conjecture of Kaplansky. They have also proved that in
characteristic $p$ sufficiently large a semisimple cosemisimple
Hopf algebra is involutory. Afterwards, using this result and a
lifting theorem, Etingof and Gelaki prove in \cite{eg} that the
antipode of a semisimple
cosemisimple Hopf algebra over any field is an involution.\\
${\;\;\;}$
Trying to generalize the above results for quasi-Hopf algebras, the first problem
which occur is: what could be an involutory quasi-Hopf algebra? We believe that
we cannot keep the same definition as in the Hopf case because, in general, $S^2$
is not a coalgebra morphism, while $id_H$ is. So one of the purposes of this paper
is to find a plausible definition for this notion. Toward this end we will use
a categorical point of view due to Majid \cite{m2}. More exactly, he has observed
that ${\rm Tr}(S^2)$, the trace of $S^2$, is an important invariant of any
finite dimensional Hopf algebra. In fact, he has shown that ${\rm Tr}(S^2)$ arises
in a very natural way as the representation-theoretic
rank of the Schr${\rm{\ddot{o}}}$dinger representation of $H$, $\un{\rm dim}(H)$, or as the
representation-theoretic rank of the canonical representation of the quantum double,
$\un{\rm dim}(D(H))$. Correlating this with the trace formula obtained by
Radford in \cite{r2} we get that
\[
\un{\rm dim}(H)=\un{\rm dim}(D(H))={\rm Tr}(S^2)=\va (r)\l (1),
\]
where $\l$ is a left integral in $H^*$ and $r$ is a right integral in $H$ such
that $\l (S(r))=1$. By the Larson-Radford-Etingof-Gelaki results we conclude that
\[
\un{\rm dim}(H)=\un{\rm dim}(D(H))=\left\{\begin{array}{ll}
0 &\hspace*{-3mm}\mbox{,~~if $H$ is not semisimple or cosemisimple}\\
{\rm dim}(H)&\hspace*{-3mm} \mbox{,~~if $H$ is both semisimple and cosemisimple.}
\end{array}\right.
\]
${\;\;\;}$ 
The aim of this paper is to generalize some of the
results presented above for quasi-Hopf algebras by computing the
representation-theoretic rank of a finite dimensional quasi-Hopf
algebra $H$ and of its quantum double $D(H)$. We hope that the
point of view presented here will open the way for solving the
remaining ones. The paper is organized as follows. In \seref{3} we
compute the Schr${\rm{\ddot{o}}}$dinger representation associated
to a finite dimensional quasi-Hopf algebra $H$. In fact, we will
transfer the associated algebra structure of $H$ within the
category of left Yetter-Drinfeld modules constructed in
\cite{bpv2, bn2} to the category of left $D(H)$-modules, through 
some monoidal isomorphisms explicitly constructed in \cite{bcp}
and \cite{bpv}. Now following \cite{m2}, in any braided rigid
monoidal category ${\cal C}$ we can compute the
representation-theoretic rank of an object $V$ of ${\cal C}$.
Considering ${\cal C}={}_{D(H)}{\cal M}^{\rm fd}$, the category of
finite dimensional left $D(H)$-modules, we will compute in
\seref{4} the representation-theoretic rank of $H$ and $D(H)$
within ${\cal C}$, $\un{\rm dim}(H)$ and $\un{\rm dim}(D(H))$,
respectively. After some technical and complicated computations we
will find that
\[
\un{\rm dim}(H)=\un{\rm dim}(D(H))={\rm Tr}\left(h\mapsto S^{-2}(S(\b)\a h\b S(\a))\right).
\]
Therefore, we call a quasi-Hopf algebra $H$ involutory 
if $h\mapsto S^{-2}(S(\b)\a h\b S(\a))=id_H$. Firstly, because, 
just as in the Hopf case, the above representation-theoretic ranks reduce to
the classical dimension of $H$, provided $H$ involutory. Secondly, because 
${\mf g}=\b S(\a)$ is invertible with ${\mf g}^{-1}=S(\b)\a$ and ${\mf g}^{-1}$ defines 
both $S^2$ as a inner automorphism of $H$ and (assuming $k$ algebraically closed of characteristic 
zero) that unique pivotal structure in \cite[Propositions 8.24 and 8.23]{eno}. 
More explicitly, if $k$ is an algebraically closed field of characteristic zero then 
${\mf g}^{-1}$ gives rise to the unique pivotal structure of ${}_H{\cal M}^{\rm fd}$ 
with respect to which the categorical dimensions of simple objects coincide with their usual dimensions.  
(Complete proofs for the above facts, examples, properties and results on 
involutory (dual) quasi-Hopf algebras can be found in \cite{bct}.)  
Furthermore, specializing the above equality for $H=H^*_{\omega}$,
the quasi-Hopf algebra considered in \cite{pv}, we obtain that
$\un{\rm dim}(D^{\omega}(H))={\rm dim}(H)$, where $D^{\omega}(H)$
is the quasi-triangular quasi-Hopf algebra constructed in
\cite{bp}, and we should stress the fact that in this particular
case both $H^*_{\omega}$ and $D^{\omega}(H)$ are involutory in the quasi-Hopf 
sense mentioned above.\\ 
${\;\;\;}$
Finally, in \seref{5} we prove a trace formula for
quasi-Hopf algebras. Specializing it for the endomorphism
$h\mapsto S^{-2}(S(\b)\a h\b S(\a))$ we get that
\[
{\rm Tr}\left(h\mapsto S^{-2}(S(\b)\a h\b S(\a))\right)=\va(r)\l (\smi (\a)\b),
\]
where $\l $ is a left cointegral in $H$ and $r$ is a right integral in $H$ such that
$\l (r)=1$. Combining the results in the last two Sections we conclude that
$\un{\rm dim}(H)=\un{\rm dim}(D(H))=\va (r)\l (\smi (\a )\b)$, so this scalar is non-zero
if and only if $H$ is both semisimple and cosemisimple.\\
${\;\;\;}$
In view of these results we believe that a semisimple cosemisimple quasi-Hopf algebra
is always involutory and therefore, in this case, $\un{\rm dim}(H)=\un{\rm dim}(D(H))=
{\rm dim}(H)$, the classical dimension of $H$. In this direction we do not know if
the techniques used in \cite{lr1,lr2,eg} can be generalized for quasi-Hopf algebras.
But without doubt it is an interesting problem which is worthwhile to study.
\section{Preliminaries}\selabel{2}
\subsection{Quasi-Hopf algebras}
We work over a commutative field $k$. All algebras, linear
spaces etc. will be over $k$; unadorned $\ot $ means $\ot_k$.
Following Drinfeld \cite{d1}, a quasi-bialgebra is
a four-tuple $(H, \Delta , \va , \Phi )$ where $H$ is
an associative algebra with unit,
$\Phi$ is an invertible element in $H\ot H\ot H$, and
$\Delta :\ H\ra H\ot H$ and $\va :\ H\ra k$ are algebra
homomorphisms satisfying the identities
\begin{eqnarray}
&&(id \ot \Delta )(\Delta (h))=
\Phi (\Delta \ot id)(\Delta (h))\Phi ^{-1},\label{q1}\\
&&(id \ot \va )(\Delta (h))=h\ot 1,
\mbox{${\;\;\;}$}
(\va \ot id)(\Delta (h))=1\ot h,\label{q2}
\end{eqnarray}
for all $h\in H$, and
$\Phi$ has to be a $3$-cocycle, in the sense that
\begin{eqnarray}
&&(1\ot \Phi)(id\ot \Delta \ot id)
(\Phi)(\Phi \ot 1)=
(id\ot id \ot \Delta )(\Phi )
(\Delta \ot id \ot
id)(\Phi ),\label{q3}\\
&&(id \ot \va \ot id)(\Phi )=1\ot 1\ot 1.\label{q4}
\end{eqnarray}
The map $\Delta $ is called the coproduct or the
comultiplication, $\va $ the counit and $\Phi $ the
reassociator. As for Hopf algebras we denote $\Delta (h)=h_1\ot h_2$,
but since $\Delta $ is only quasi-coassociative we adopt the
further convention (summation understood):
$$
(\Delta \ot id)(\Delta (h))=h_{(1, 1)}\ot h_{(1, 2)}\ot h_2,
\mbox{${\;\;\;}$}
(id\ot \Delta )(\Delta (h))=h_1\ot h_{(2, 1)}\ot h_{(2,2)},
$$
for all $h\in H$. We will
denote the tensor components of $\Phi$
by capital letters, and the ones of
$\Phi^{-1}$ by small letters,
namely
\begin{eqnarray*}
&&\Phi=X^1\ot X^2\ot X^3=T^1\ot T^2\ot T^3=
V^1\ot V^2\ot V^3=\cdots\\
&&\Phi^{-1}=x^1\ot x^2\ot x^3=t^1\ot t^2\ot t^3=
v^1\ot v^2\ot v^3=\cdots
\end{eqnarray*}
$H$ is called a quasi-Hopf
algebra if, moreover, there exists an
anti-morphism $S$ of the algebra
$H$ and elements $\a , \b \in
H$ such that, for all $h\in H$, we
have:
\begin{eqnarray}
&&S(h_1)\a h_2=\va(h)\a
\mbox{${\;\;\;}$ and ${\;\;\;}$}
h_1\b S(h_2)=\va (h)\b,\label{q5}\\
&&X^1\b S(X^2)\a X^3=1
\mbox{${\;\;\;}$ and${\;\;\;}$}
S(x^1)\a x^2\b S(x^3)=1.\label{q6}
\end{eqnarray}
Our definition of a quasi-Hopf algebra is different from the
one given by Drinfeld \cite{d1} in the sense that we do not
require the antipode to be bijective. Nevertheless, in the finite
dimensional or quasi-triangular case this condition can be deleted because
it follows from the other axioms, see \cite{bc1} and \cite{bn3}.\\
${\;\;\;}$
Together with a quasi-Hopf algebra
$H=(H, \Delta , \va , \Phi , S, \a , \b )$ we also have $H^{\rm op}$
and $H^{\rm cop}$ as quasi-Hopf algebras, where "op" means opposite
multiplication and "cop" means opposite comultiplication. The quasi-Hopf
structures are obtained by putting $\Phi_{\rm op}=\Phi^{-1}$,
$\Phi_{\rm cop}=(\Phi ^{-1})^{321}$, $S_{\rm op}=S_{\rm cop}=S^{-1}$,
$\a _{\rm op}=\smi (\b )$, $\b _{\rm op}=\smi (\a )$, $\a _{\rm cop}=\smi (\a )$
and $\b _{\rm cop}=\smi (\b )$.\\
${\;\;\;}$
The axioms for a quasi-Hopf algebra imply that
$\va \circ S=\va $ and $\va (\a )\va (\b )=1$,
so, by rescaling $\a $ and $\b $, we may assume without loss of generality
that $\va (\a )=\va (\b )=1$. The identities
(\ref{q2}), (\ref{q3}) and (\ref{q4}) also imply that
\begin{equation}\label{q7}
(\va \ot id\ot id)(\Phi )=
(id \ot id\ot \va )(\Phi )=1\ot 1\ot 1.
\end{equation}
${\;\;\;}$
It is well-known that the antipode of a Hopf algebra
is an anti-coalgebra morphism. For a quasi-Hopf algebra, we have
the following statement: there exists an invertible element
$f\in H\ot H$ such that $(\va \ot id)(f)=(id\ot \va)(f)=1$ and
\begin{equation} \label{ca}
f\Delta (S(h))f^{-1}=(S\ot S)(\Delta ^{op}(h))
\mbox{,${\;\;\;}$for all $h\in H$,}
\end{equation}
where $\Delta ^{op}(h)=h_2\ot h_1$. $f$ can be computed explicitly. First set
\begin{eqnarray*}
A^1\ot A^2\ot A^3\ot A^4&=&
(\Phi \ot 1) (\Delta \ot id\ot id)(\Phi
^{-1}),\\
B^1\ot B^2\ot B^3\ot B^4&=&
(\Delta \ot id\ot id)(\Phi )(\Phi ^{-1}\ot 1)
\end{eqnarray*}
and then define $\gamma, \delta\in H\ot H$ by
\begin{equation} \label{gd}%
\gamma =S(A^2)\a A^3\ot S(A^1)\a A^4~~{\rm and}~~
\delta =B^1\b S(B^4)\ot B^2\b S(B^3).
\end{equation}
$f$ and $f^{-1}$ are then given by
the formulas
\begin{eqnarray}
f&=&(S\ot S)(\Delta ^{op}(x^1)) \gamma
\Delta (x^2\b S(x^3)),\label{f}\\
f^{-1}&=&\Delta (S(x^1)\a x^2)
\delta (S\ot S)(\Delta ^{op}(x^3)).\label{g}
\end{eqnarray}
Moreover, $f=f^1\ot f^2$ and $f^{-1}=g^1\ot g^2$ satisfy the following relations:
\begin{eqnarray}
&&\hspace*{5mm}
f\Delta (\a )=\gamma ,~~
\Delta (\b )f^{-1}=\delta ,\label{gdf}\\
&&\hspace*{5mm}
(1\ot f)(id \ot \Delta )(f) \Phi (\Delta \ot
id)(f^{-1})(f^{-1}\ot 1)
=S(X^3)\ot S(X^2)\ot S(X^1),\label{pf}\\
&&\hspace*{5mm}
f^1\b S(f^2)=S(\a ),~~
g^1S(g^2\a )=\b ,~~
S(\b f^1)f^2=\a .\label{fgab}
\end{eqnarray}
${\;\;\;}$
In a Hopf algebra $H$, we obviously have the identity
\[
h_1\ot h_2S(h_3)=h\ot 1,~{\rm for~all~}h\in H.
\]
We will need the generalization of this formula to
quasi-Hopf algebras. Following \cite{hn1, hn2}, we define
\begin{eqnarray}
&&p_R=p^1\ot p^2=x^1\ot x^2\b S(x^3),\hspace*{5mm}
q_R=q^1\ot q^2=X^1\ot S^{-1}(\a X^3)X^2,\label{qr}\\
&&p_L=\tpla \ot \tplb=X^2\smi (X^1\b )\ot X^3,\hspace*{5mm}
q_L=\tqla \ot \tqlb=S(x^1)\a x^2\ot x^3.\label{ql}
\end{eqnarray}
For all $h\in H$, we then have
\begin{eqnarray}
\Delta (h_1)p_R(1\ot S(h_2))&=&p_R(h\ot 1)\label{qr1}\\
(S(h_1)\ot 1)q_L\Delta (h_2)&=&(1\ot h)q_L.\label{ql1a}
\end{eqnarray}
Furthermore, the following relations hold
\begin{eqnarray}
&&(1\ot S^{-1}(p^2))q_R\Delta (p^1)=1\ot 1\label{pqra}\\
&&\Delta (q^1)p_R(1\ot S(q^2))=1\ot 1\label{pqr}\\
&&(S(\tpla)\ot 1)q_L\Delta (\tplb)=1\ot 1\label{pql}\\
&&\Phi (\Delta \ot id)(p_R)(p_R\ot id)\nonumber\\
&&\hspace*{1cm}
=(id\ot \Delta )(\Delta (x^1)p_R)(1\ot f^{-1})(1\ot S(x^3)\ot
S(x^2))\label{pr}\\
&&(q_R\ot 1)(\Delta \ot id)(q_R)\Phi ^{-1}\nonumber\\
&&\hspace*{1cm}
=(1\ot S^{-1}(f^2X^3)\ot S^{-1}(f^1X^2))
(id \ot \Delta )(q_R\Delta (X^1)),\label{qr2}\\
&&(1\ot q_L)(id \ot \Delta )(q_L)\Phi\nonumber\\
&&\hspace*{1cm}
=(S(x^2)\ot S(x^1)\ot 1)(f\ot 1)(\Delta \ot id)(q_L\Delta (x^3)).\label{ql2}
\end{eqnarray}
\subsection{Quasi-triangular quasi-Hopf algebras and the quantum double}
Recall that a quasi-Hopf algebra $H$ is quasi-triangular if
there exists an element $R\in H\ot H$ such that
\begin{eqnarray}
(\Delta \ot id)(R)&=&\Phi _{312}R_{13}\Phi ^{-1}_{132}R_{23}\Phi
,\label{qt1}\\
(id \ot \Delta )(R)&=&\Phi ^{-1}_{231}R_{13}\Phi _{213}R_{12}\Phi
^{-1},\label{qt2}\\
\Delta ^{\rm op}(h)R&=&R\Delta (h),~{\rm for~all~}h\in H,\label{qt3}\\
(\va \ot id)(R)&=&(id\ot \va)(R)=1.\label{qt4}
\end{eqnarray}
Here we use the following notation.
If $\sigma $ is a permutation of $\{1, 2, 3\}$, we set $\Phi _{\sigma (1)
\sigma (2)\sigma (3)}=X^{\sigma ^{-1}(1)}\ot
X^{\sigma ^{-1}(2)}\ot X^{\sigma ^{-1}(3)}$, and
$R_{ij}$ means $R$ acting non-trivially
in the $i^{th}$ and $j^{th}$ positions of $H\ot H\ot H$.\\
${\;\;\;}$
In {\cite{bn3}} it is shown that $R$ is invertible, and that
the element
\begin{equation} \label{elmu}
u=S(R^2p^2)\a R^1p^1
\end{equation}
(with $p_R=p^1\ot p^2$ defined as in (\ref{qr})) is
invertible in $H$ and satisfies for all $h\in H$ the following relation
\begin{equation} \label{sqina}
S^2(h)=uhu^{-1}.
\end{equation}
${\;\;\;}$
As in the Hopf algebra theory the most important example of quasi-triangular
quasi-Hopf algebra is produced by the double construction.\\
${\;\;\;}$
From \cite{hn2, bc}, we recall the definition of the quantum double
$D(H)$ of a finite dimensional quasi-Hopf algebra $H$.
Let $\{e_i\}_{i=\ov {1, n}}$ be a basis of $H$, and $\{e^i\}_{i=\ov
{1, n}}$ the corresponding dual basis of $H^*$. We can easily see that $H^*$,
the linear dual of $H$, is not a quasi-Hopf algebra.
But $H^*$ has a dual structure coming from the initial structure of $H$.
So $H^*$ is a coassociative coalgebra, with comultiplication
\[
\widehat {\Delta }(\v )=\v _1\ot \v _2=
\sum \limits _{i, j=1}^n\v (e_ie_j)e^i\ot e^j,%
\]
or, equivalently, %
\[
\widehat {\Delta }(\v )=\v _1\ot \v_2 \Leftrightarrow
\v (hh')=\v _1(h)\v _2(h'), %
\mbox{${\;\;\;}$$\forall ~~h, h'\in H$.}
\]
$H^*$ is also an $H$-bimodule, by
\[
\<h\rh \v , h'\>=\v (h'h),
\mbox{${\;\;\;}$}
\<\v \lh h, h'\>=\v (hh').
\]
The convolution is a multiplication on $H^*$; it is not
associative, but only quasi-associative:
\[
[\v \psi]\xi=(X^1\rh \v \lh x^1)[(X^2\rh \psi \lh x^2)
(X^3\rh \xi \lh x^3)],~~\forall ~~\v , \psi , \xi \in H^*.
\]
We also introduce $\ov {S}:\ H^*\ra H^*$ as the coalgebra
antimorphism dual to $S$, this means
$\<\ov {S}(\v), h\>=$ $\<\v, S(h)\>$, for all $\v \in H^*$ and $h\in H$.\\
${\;\;\;}$
Now consider $\Omega \in H^{\ot 5}$ given by
\begin{eqnarray}
&&\hspace*{-5mm} \Omega =\Omega^1\ot \Omega^2\ot \Omega^3\ot
\Omega^4\ot \Omega^5 \nonumber \\%
&&\hspace*{5mm}
=X^1_{(1, 1)}y^1x^1\ot X^1_{(1, 2)}y^2x^2_1\ot
X^1_2y^3x^2_2\ot \smi (f^1X^2x^3)\ot \smi (f^2X^3), \label{O}
\end{eqnarray}
where $f\in H\ot H$ is the element defined in (\ref{f}). We define
the quantum double $D(H)=H^*\bowtie H$ as follows: as a $k$-linear
space, $D(H)$ equals $H^*\ot H$, and the multiplication is given
by
\begin{eqnarray}
&&\hspace*{-1cm}(\v \bowtie h)(\psi \bowtie h')\nonumber\\
&=&[(\Omega ^1\rh \v \lh \Omega ^5)(\Omega ^2\rh \psi _2\lh \Omega
^4)]\bowtie \Omega ^3[(\ov {S}^{-1}(\psi _1)\rh h)\lh \psi
_3]h'\nonumber\\
&=&[(\Omega ^1\rh \v \lh \Omega ^5)(\Omega
^2h_{(1, 1)}\rh \psi \lh \smi (h_2)\Omega ^4)]\bowtie \Omega
^3h_{(1, 2)}h'.\label{mdd}
\end{eqnarray}

From \cite{hn1,hn2} we have that $D(H)$ is an
associative algebra with unit $\va \bowtie 1$, and $H$ is a unital
subalgebra via the morphism $i_D:\ H\ra D(H)$, $i_D(h)=\va \bowtie
h$. Moreover, $D(H)$ is a quasi-triangular quasi-Hopf algebra
with the following structure:
\begin{eqnarray}
&&\hspace*{-5mm}\Delta _D(\v \bowtie h)=
(\va \bowtie X^1Y^1)
(p^1_1x^1\rh \v _2\lh Y^2\smi (p^2)\bowtie p^1_2x^2h_1)\nonumber\\
&&\hspace*{10mm}\ot (X^2_1\rh \v _1\lh \smi (X^3)\bowtie X^2_2Y^3x^3h_2)
\label{cddf}\\
&&\hspace*{-5mm}
\va _D(\v \bowtie h)=\va (h)\v (\smi (\a ))\label{coddf}\\
&&\hspace*{-5mm}
\Phi _D=(i_D\ot i_D\ot i_D)(\Phi)\label{pdd}\\
&&\hspace*{-5mm}
S_D(\v \bowtie h)=(\va \bowtie S(h)f^1)(p^1_1U^1\rh \ov
{S}^{-1}(\v )\lh f^2\smi (p^2)\bowtie p^1_2U^2)\label{anddf}\\
&&\hspace*{-5mm}
\a _D=\va \bowtie \a ,~~\b _D=\va \bowtie \b \label{abdd}\\
&&\hspace*{-5mm}
R_D=\sum \limits _{i=1}^n
(\va \bowtie \smi (p^2)e_ip^1_1)\ot (e^i\bowtie p^1_2).\label{rdd}
\end{eqnarray}
Here $p_R=p^1\ot p^2$ and $f=f^1\ot f^2$ are the elements defined by
(\ref{qr}) and (\ref{f}), respectively, and $U=U^1\ot U^2\in H\ot H$
is the following element
\begin{equation}\label{U}
U=U^1\ot U^2=g^1S(q^2)\ot g^2S(q^1),
\end{equation}
where $f^{-1}=g^1\ot g^2$ and $q_R=q^1\ot q^2$ are the elements
defined by (\ref{g}) and (\ref{qr}), respectively.

\subsection{The center construction and the Yetter-Drinfeld modules}
If $H$ is a quasi-bialgebra then the category of left $H$-modules,
denoted by ${}_H{\cal M}$, is a monoidal category and, moreover,
if $H$ is quasi-triangular then ${}_H{\cal M}$ is braided (the reader is
invited to consult \cite[XI.4]{k} or \cite[IX.1]{maj} for the complete
definition of a monoidal or (pre) braided category, and also for the notion
of a monoidal, respectively (pre) braided, functor between them).
The tensor product $\ot $ is given via $\Delta $, for $U,V,W\in {}_H\Mm$
the associativity constraint on ${}_H\Mm$ is given by
\begin{equation}\eqlabel{1.3.1}
a_{U,V,W}((u\ot v)\ot w)=X^1\cd u\ot (X^2\cd v\ot X^3\cd w),
\end{equation}
the unit is $k$ as a trivial $H$-module and the left and right unit
constraints are the usual ones. When $H$ is quasi-triangular we have the following
braiding $c$ on ${}_H\Mm$:
\begin{equation}\label{br}
c_{U, V}(u\ot v)=R^2\cd v\ot R^1\cd u.
\end{equation}
${\;\;\;}$
To any monoidal category ${\cal C}$ we can associate two
(pre) braided monoidal categories, namely the (weak) left
and right centers (${\cal W}_{l/r}({\cal C})$)
${\cal Z}_{l/r}({\cal C})$ of ${\cal C}$.
For the (weak) left center construction the reader is invited
to consult \cite{m1}, for the right (weak) center construction
\cite[XIII.4]{k}, and for the connection between them \cite{bcp},
respectively.\\
${\;\;\;}$
Since for a quasi-bialgebra $H$ the category ${}_H{\cal M}$ is monoidal
it makes sense to consider ${\cal W}_l({}_H{\cal M})$ or ${\cal W}_r({}_H{\cal M})$.
In \cite{m1} Majid computed the left weak center ${\cal W}_l({}_H{\cal M})$.
The objects are identified with the so called left Yetter-Drinfeld modules,
i.e. left $H$-modules $M$ (denote the action by $h\ot m\mapsto h\cd m$)
together with a $k$-linear map $\l _M: M\ra H\ot M$,
$\l _M(m):=m_{\<-1\>}\ot m_{\<0\>}$, such that $\va (m_{\<-1\>})m_{\<0\>}=m$ and
for all $h\in H$ and
$m\in M$ the following relations hold:
\begin{eqnarray}
&&X^1m_{\<-1\>}\ot (X^2\cd m_{\<0\>})_{\<-1\>}X^3
\ot (X^2\cd m_{\<0\>})_{\<0\>}\nonumber\\
&&\hspace*{1cm}=X^1(Y^1\cd m)_{\<-1\>_1}Y^2\ot X^2(Y^1\cd
m)_{\<-1\>_2}Y^3
\ot X^3\cd (Y^1\cd m)_{\<0\>},\label{y1}\\
&&h_1m_{\<-1\>}\ot h_2\cd m_{\<0\>}=(h_1\cd
m)_{\<-1\>}h_2\ot (h_1\cd m)_{\<0\>}.\label{y3}
\end{eqnarray}
The category of left Yetter-Drinfeld modules and $k$-linear
maps that preserve the $H$-action and $H$-coaction is denoted by
$\yd$. \\
${\;\;\;}$
The prebraided monoidal structure on ${\cal W}_l({}_H\Mm)$ induces a prebraided
monoidal structure on $\yd$. This structure is such that the forgetful functor
$\yd\to {}_H{\mathcal M}$ is monoidal, and the coaction on the tensor product
$M\ot N$ of two left Yetter-Drinfeld modules $M$ and $N$ is given by
\begin{eqnarray}
&&\hspace*{-2cm}
\lambda _{M\ot N}(m\ot n)=X^1(x^1Y^1\cd m)_
{\<-1\>} x^2(Y^2\cd n)_{\<-1\>}Y^3\nonumber\\
&&\hspace*{1.5cm}
\ot X^2\cd (x^1Y^1\cd
m)_{\<0\>}\ot X^3x^3\cd (Y^2\cd n)_{\<0\>}.\label{y4}
\end{eqnarray}
For any $M, N\in \yd$ the braiding $c_{M, N}: M\ot N\ra N\ot M$ is given by
\begin{equation}\label{y5}
c_{M, N}(m\ot n)=m_{\<-1\>}\cd n\ot m_{\<0\>},
\end{equation}
for all $m\in M$ and $n\in N$. Moreover, if $H$ is a quasi-Hopf algebra then
$c_{M, N}$ is invertible (see \cite{bn2}) and therefore $\yd$ is a braided category.\\
${\;\;\;}$
We notice that the right weak center ${\cal W}_r({}_H{\cal M})$
was computed in \cite{bcp}: it is isomorphic to the category of left-right
Yetter-Drinfeld modules (see the definition below).

\section{The Schr${\rm{\ddot{o}}}$dinger representation}\selabel{3}
\setcounter{equation}{0}
${\;\;\;}$
Let $H$ be a finite dimensional Hopf algebra. It is well known that $H$ is a
left $D(H)$-module algebra via the action ($\v \in H^*$, $h, h'\in H$):
\[
(\v\bowtie h)\bullet h'=\<\v , \smi \left((h\tr h')_1\right)\>(h\tr h')_2.
\]
Here and also in the rest of the paper, $h\tr h':=h_1h'S(h_2)$, for all $h, h'\in H$.\\
${\;\;\;}$
The aim of this section is to compute a similar structure for a finite
dimensional quasi-Hopf algebra $H$. This fact is absolutely necessary in order
to compute the representation-theoretic rank (or quantum dimension) of $H$ within the
braided category of left $D(H)$-modules. Toward this end, we will use the
following three results:
\begin{itemize}
\item[1)] To any quasi-Hopf algebra $H$ we can associate an algebra, denoted by
$H_0$, in the category of left Yetter-Drinfeld modules $\yd$, cf. \cite{bn2}.
\end{itemize}
More precisely, we denote by $H_0$ the $k$-vector space $H$ with the new multiplication
$\circ$ defined by
\begin{equation}\label{sma}
h\circ h'=X^1hS(x^1X^2)\a x^2X^3_1h'S(x^3X^3_2),
\end{equation}
for all $h, h'\in H$. From \cite{bpv2} we know that $H_0$ is a left
$H$-module algebra, this means an algebra in ${}_H{\cal M}$.
The unit of $H_0$ is $\b$ and $H_0$ is an object of
${}_H{\cal M}$ via the left adjoint action $\tr$, i.e. for all
$h, h'\in H$,
\begin{equation}\label{lya2}
h\tr h'=h_1h'S(h_2).
\end{equation}
Moreover, $H_0$ becomes an algebra in $\yd$ with the additional structure
$\l _{H_0}:\ H_0\ra H\ot H_0$, given by
\begin{equation}\label{lyda3}
\l _{H_0}(h)=h_{\<-1\>}\ot h_{\<0\>}:=
X^1Y^1_1h_1g^1S(q^2Y^2_2)Y^3\ot X^2Y^1_2h_2g^2S(X^3q^1Y^2_1),
\end{equation}
for all $h\in H$, where $q_R=q^1\ot q^2$ is the element defined by (\ref{qr}).
\begin{itemize}
\item[2)] There is a braided isomorphism between
$\yd$ and ${{}_H{\cal YD}^H}^{\rm in}$, cf. \cite{bcp}.
\end{itemize}
First of all recall that the category of left-right Yetter-Drinfeld
modules over a quasi-bialgebra $H$,
denoted by ${}_H{\cal YD}^H$, has as objects left $H$-modules
$M$ (denote the action by $h\ot m\mapsto h\cd m$) for which $H$ coacts on
the right (denote the right $H$-coaction by $M\ni m\mapsto
m_{(0)}\ot m_{(1)}\in M\ot H$) such that $\varepsilon (m_{(1)})m_{(0)}=m$
and for all $m\in M$ and $h\in H$ the
following relations hold
\begin{eqnarray}
&&\hspace*{-1.5cm}
(x^2\cdot m_{(0)})_{(0)}\otimes (x^2\cdot m_{(0)})_{(1)}x^1
\otimes x^3m_{(1)}\nonumber \\
&=&x^1\cdot (y^3\cdot m)_{(0)}\otimes x^2(y^3\cdot m)_{(1)_1}y^1
\otimes x^3(y^3\cdot m)_{(1)_2}y^2,\label{lry1}\\
&&\hspace*{-1.5cm}
h_1\cdot m_{(0)}\otimes h_2m_{(1)}=(h_2\cdot m)_{(0)}
\otimes (h_2\cdot m)_{(1)}h_1\label{lry3}.
\end{eqnarray}
The morphisms are left $H$-linear, right $H$-colinear maps. \\
${\;\;\;}$
Since ${}_H{\cal YD}^H$ can be identified with the right weak center
of ${}_H{\cal M}$ (see \cite{bcp}) we find that ${}_H{\cal YD}^H$ has the following
prebraided structure: the right $H$-coaction on the tensor product $M{\ot}N$ of
$M, N\in {}_H{\cal YD}^H$ is the following:
\begin{eqnarray}
&&\hspace*{-2cm}
{\r }_{M{\ot}N}(m{\ot} n)=x^1X^1\cd (y^2\cd m)_{(0)}
{\ot} x^2\cd (X^3y^3\cd n)_{(0)}\nonumber\\
&&\hspace*{2cm}
\ot x^3(X^3y^3\cd n)_{(1)}X^2(y^2\cd m)_{(1)}y^1,\label{slrms2}
\end{eqnarray}
for all $m\in M$, $n\in N$, and the functor forgetting the
$H$-coaction is monoidal, so
\begin{equation}
h\cd (m{\ot} n)=h_1\cd m{\ot} h_2\cd n.\label{slrms1}
\end{equation}
The braiding ${\mf{c}}$ on ${}_H{\cal YD}^H$ is defined by
${\mf{c}}_{M, N}: M{\ot}N\ra N{\ot}M$,
\begin{equation}\label{slrbs}
{\mf{c}}_{M, N}(m{\ot}n)=n_{(0)}{\ot}n_{(1)}\cd m,
\end{equation}
for $m\in M$ and $n\in N$. Furthermore, if $H$ is a quasi-Hopf algebra the braiding
${\mf{c}}$ is invertible. The inverse braiding is given by
\begin{equation}\label{slribs}
{\mf{c}}_{M, N}^{-1}(n{\ot} m)=q^1_1x^1
S(q^2x^3(\tplb \cd n)_{(1)}\tpla )\cd m{\ot}q^1_2x^2\cd (\tplb \cd n)_{(0)},
\end{equation}
where $q_R=q^1\ot q^2$ and $p_L=\tpla \ot \tplb$ are the elements defined in
(\ref{qr}) and (\ref{ql}), respectively.
Finally, for a quasi-Hopf algebra $H$, ${}_H{\cal YD}^H$ will be our notation
for the category of left-right Yetter-Drinfeld modules endowed with the braided structure
given by (\ref{slrms2}-\ref{slrbs}), and
${{}_H{\cal YD}^H}^{\rm in}$ will be our notation for
the category ${}_H{\cal YD}^H$ with monoidal structure (\ref{slrms2}-\ref{slrms1})
and the mirror reversed braiding
$\tilde{{\mf c}}_{M,N}={\mf{c}}^{-1}_{N, M}$.\\
${\;\;\;}$
Now, by \cite{bcp} there is a  monoidal isomorphism between
$\yd$ and ${}_H{\cal YD}^H$ produced
by the following functor $\mf{F}$. If $M\in \yd$ then $\mf{F}(M)=M$
as left $H$-modules and with the right $H$-coaction defined by
\begin{equation}\label{funct}
\r _{F(M)}(m)=\tilde{q}^2_1X^2\cd (p^1\cd m)_{\<0\>}\ot \tilde{q}^2_2X^3
\smi (\tqla X^1(p^1\cd m)_{\<-1\>}p^2),
\end{equation}
for all $m\in M$. The functor $\mf{F}$ acts as identity on morphisms. Moreover,
$\mf{F}$ provides a braided isomorphism between $\yd$ and
${{}_H{\cal YD}^H}^{\rm in}$, see \cite{bcp} for more details.
\begin{itemize}
\item[3)] The category ${}_H{\cal YD}^H$ is braided isomorphic to
${}_{D(H)}{\cal M}$.
\end{itemize}
Indeed, from \cite{hn2, bpv} we know that the above categories are isomorphic.
The isomorphism is the following. To any left-right Yetter-Drinfeld module $M$
we can associate a left $D(H)$-module structure given by
\begin{equation}\label{drs}
(\v \bowtie h)\act m=\< \v , q^2(h\cd m)_{(1)}\> q^1\cd (h\cd m)_{(0)},
\end{equation}
for all $\v \in H^*$, $h\in H$ and $m\in M$, where, as usual,
$q_R=q^1\ot q^2$ is the element defined in (\ref{qr}). Moreover, a morphism
between two left-right Yetter-Drinfeld modules becomes in this way a morphism
between two left $D(H)$-modules, so we have a well defined functor
${\cal F}: {}_H{\cal YD}^H\ra {}_{D(H)}{\cal M}$. If $H$ is a finite dimensional
quasi-Hopf algebra then ${\cal F}$ is an isomorphism (for the explicit description
of the inverse of ${\cal F}$ see \cite{bpv}). Also, it is not hard to see that
the functor ${\cal F}$ is monoidal; the functorial isomorphism
$\Psi _{M, N}: {\cal F}(M)\ot {\cal F}(N)\ra {\cal F}(M\ot N)$ is the identity
morphism. Moreover, the next result asserts that it is a braided isomorphism.

\begin{proposition}\label{pr3.1}
Let $H$ be a finite dimensional quasi-Hopf algebra. Then the categories
${}_H{\cal YD}^H$ and ${}_{D(H)}{\cal M}$ are braided isomorphic.
\end{proposition}
\begin{proof}
By the previous comments, we only have to check that the functor
${\cal F}$ defined above is braided, this means that for any two left-right
Yetter-Drinfeld modules $M$ and $N$ we have
\[
{\cal F}(\mf{c}_{M, N})\circ
\Psi _{{\cal F}(M), {\cal F}(N)}=\Psi _{{\cal F}(N), {\cal F}(M)}
\circ c_{{\cal F}(M), {\cal F}(N)}.
\]
Indeed, for all $m\in M$ and $n\in N$ we compute
\begin{eqnarray*}
&&\hspace*{-3cm}\Psi _{{\cal F}(N), {\cal F}(M)}
\circ c_{{\cal F}(M), {\cal F}(N)}(m\ot n)\\
&{{\rm (\ref{br}, \ref{rdd})}\atop =}&
\sum \limits _{i=1}^n (e^i\bowtie p^1_2)\act n\ot (\va \bowtie \smi (p^2)e_ip^1_1)\act m\\
&{{\rm (\ref{drs})}\atop =}&\<e^i, q^2(p^1_2\cd n)_{(1)}\>
q^1\cd (p^1_2\cd n)_{(0)}\ot \smi (p^2)e_ip^1_1\cd m\\
&=&q^1\cd (p^1_2\cd n)_{(0)}\ot \smi (p^2)q^2(p^1_2\cd n)_{(1)}p^1_1\cd m\\
&{{\rm (\ref{lry3})}\atop =}&
q^1p^1_1\cd n_{(0)}\ot \smi (p^2)q^2p^1_2n_{(1)}\cd m\\
&{{\rm (\ref{pqra}, \ref{slrbs})}\atop =}&
n_{(0)}\ot n_{(1)}\cd m=
{\cal F}(\mf{c}_{M, N})\circ
\Psi _{{\cal F}(M), {\cal F}(N)}(m\ot n),
\end{eqnarray*}
as needed, so the proof is complete.
\end{proof}
${\;\;\;}$
Using these braided isomorphisms we will transfer the algebra structure of $H_0$
in $\yd$ to ${}_{D(H)}{\cal M}$. In this way we will associate to any finite
dimensional quasi-Hopf algebra $H$ a left $D(H)$-module algebra structure.
As in the classical Hopf algebra case, the obtained representation will be called
the Schr${\rm{\ddot{o}}}$dinger representation.\\
${\;\;\;}$
First we shall compute the algebra structure of $H_0$ in ${}_H{\cal YD}^H$, and then
its left $D(H)$-module algebra structure.

\begin{proposition}\label{pr3.2}
Let $H$ be a quasi-Hopf algebra. Then $H_0$ is an algebra in the monoidal
category ${}_H{\cal YD}^H$ with the left $H$-module structure defined in
(\ref{lya2}) and with the right $H$-coaction $\r _{H_0}: H_0\ra H_0\ot H$ given for
all $h\in H$ by
\begin{equation}\label{aclryd}
\r _{H_0}(h)=h_{(0)}\ot h_{(1)}=x^1\tqlb y^2_2h_2g^2S(x^2y^3_1)\ot x^3y^3_2
\smi (\tqla y^2_1h_1g^1)y^1,
\end{equation}
where $q_L=\tqla \ot \tqlb$ and $f^{-1}=g^1\ot g^2$ are the elements defined in
(\ref{ql}) and (\ref{g}), respectively. Moreover, $H_0$ is a left
$D(H)$-module algebra via the action
\begin{eqnarray}
&&\hspace*{-1cm}
(\v \bowtie h)\act h'=
\< \v , q^2x^3y^3_2\smi (\tqla y^2_1(h\tr h')_1g^1)y^1\>\nonumber\\
&&\hspace*{3cm}
q^1_1x^1\tqlb y^2_2(h\tr h')_2g^2S(q^1_2x^2y^3_1),\label{aadd}
\end{eqnarray}
for all $\v \in H^*$ and $h, h'\in H$, where $q_R=q^1\ot q^2$ is the element
defined in (\ref{qr}).
\end{proposition}
\begin{proof}
Since the functor $\mf{F}$ described in (\ref{funct}) is monoidal it carries
algebras to algebras. Moreover, the isomorphisms $\Psi _{M, N}:
\mf{F}(M)\ot  \mf{F}(N)\ra \mf{F}(M\ot N)$, $M, N\in \yd$,
which define the monoidal structure of the functor $\mf{F}$ are trivial,
so if $A$ is an algebra in $\yd$ then $\mf{F}(A)$ is an algebra in
${}_H{\cal YD}^H$ with the same multiplication and unit. Now, $\mf{F}$ acts as identity
on objects at the level of actions. Thus $\mf{F}({H_0})=H_0$ as left $H$-module algebras,
so we only have to show that that corresponding right $H$-action on $H_0$ through the functor
$\mf{F}$ is the one claimed in (\ref{aclryd}). For this we need the following relations
\begin{eqnarray}
&&X^1p^1_1\ot X^2p^1_2\ot X^3p^2=x^1\ot x^2_1p^1\ot x^2_2p^2S(x^3),\label{f1}\\
&&\tqla X^1\ot \tilde{q}^2_1X^2\ot \tilde{q}^2_2X^3=
S(x^1)\tqla x^2_1\ot \tqlb x^2_2\ot x^3,\label{f2}\\
&&f^{-1}=\Delta (S(p^1))U(p^2\ot 1),\label{f3}\\
&&\Delta (S(h_1))U(h_2\ot 1)=U(1\ot S(h)),~~\forall ~~h\in H. \label{f4}
\end{eqnarray}
Indeed, (\ref{f1}) and (\ref{f2}) follow easily from (\ref{q3}), (\ref{q5}) and
from the definitions of $p_R$ and $q_L$, respectively. The relation (\ref{f3}) is an
immediate consequence of (\ref{ca}) and (\ref{pqra}), and the formula in (\ref{f4})
can be found in \cite{hn3}.\\
${\;\;\;}$
Finally, by (\ref{ca}) and (\ref{U}) we get the following second formula for the
left $H$-coaction on $H_0$ defined in (\ref{lyda3})
\[
\l _{H_0}(h)=h_{\<-1\>}\ot h_{\<0\>}=(X^1\ot X^2)\Delta (Y^1hS(Y^2))U(Y^3\ot S(X^3)).
\]
${\;\;\;}$
Now, for any $h\in H$ we calculate
\begin{eqnarray*}
\r _{H_0}(h)
&{{\rm (\ref{funct})}\atop =}&
\tilde{q}^2_1Z^2\tr (p^1\tr h)_{\<0\>}\ot \tilde{q}^2_2Z^3
\smi (\tqla Z^1(p^1\tr h)_{\<-1\>}p^2)\\
&=&\tilde{q}^2_1Z^2\tr [X^2\left(Y^1(p^1\tr h)S(Y^2)\right)_2U^2S(X^3)]\\
&&\hspace*{1cm}
\ot \tilde{q}^2_2Z^3\smi (\tqla Z^1X^1\left(Y^1(p^1\tr h)S(Y^2)\right)_1U^1Y^3p^2)\\
&{{\rm (\ref{lya2}, \ref{f1})}\atop =}&
\tilde{q}^2_1Z^2\tr[X^2x^1_2h_2S(x^2_1p^1)_2U^2S(X^3)]\\
&&\hspace*{1cm}
\ot \tilde{q}^2_2Z^3x^3\smi (\tqla Z^1X^1x^1_1h_1S(x^2_1p^1)_1U^1x^2_2p^2)\\
&{{\rm (\ref{lya2}, \ref{f4})}\atop =}&
\tilde{q}^2_1\tr [Z^2_1X^2x^1_2h_2S(p^1)_2U^2S(Z^2_2X^3x^2)]\\
&&\hspace*{1cm}
\ot \tilde{q}^2_2Z^3x^3\smi (\tqla Z^1X^1x^1_1h_1S(p^1)_1U^1p^2)\\
&{{\rm (\ref{q3}, \ref{f3}, \ref{lya2})}\atop =}&
\tilde{q}^2_{(1, 1)}x^1X^2h_2g^2S(\tilde{q}^2_{(1, 2)}x^2X^3_1)
\ot \tilde{q}^2_2x^3X^3_2\smi (\tqla X^1h_1g^1)\\
&{{\rm (\ref{q1}, \ref{f2})}\atop =}&
x^1\tqlb y^2_2h_2g^2S(x^2y^3_1)\ot x^3y^3_2\smi (\tqla y^2_1h_1g^1)y^1,
\end{eqnarray*}
as needed. The last assertion is a consequence of (\ref{drs}) and (\ref{aclryd}),
the details are left to the reader.
\end{proof}

\begin{remark}\label{re3.3}
Let $H$ be a quasi-triangular quasi-Hopf algebra. Under this condition it was proved in
\cite{bn2} that $H_0$ is a braided Hopf algebra in $\yd$. Using the functor $\mf{F}$
described above we obtain that $H_0$ has also a braided Hopf algebra structure
in ${{}_H{\cal YD}^H}^{\rm in}$; note that the left
$H$-coaction of $H_0$ in $\yd$ (viewed as a braided Hopf algebra)
is different from the coaction defined in (\ref{lya2}), so the braided Hopf algebra
structure of $H_0$ within ${{}_H{\cal YD}^H}^{\rm in}$ is not induced by the algebra
structure of $H_0$ obtained in Proposition \ref{pr3.2}. Furthermore,
if we want to associate to $H$ a braided Hopf algebra in
${}_H{\cal YD}^H$ (and therefore in ${}_{D(H)}{\cal M}$ when $H$ is
finite dimensional), is sufficient to consider $H_0^{\rm op}$ (or $H_0^{\rm cop}$),
the opposite (the coopposite, respectively) braided Hopf algebra associated to $H_0$.
We leave the verification of the details to the reader.
\end{remark}

\section{The representation-theoretic rank}\selabel{4}
\setcounter{equation}{0}
${\;\;\;}$
Let ${\cal C}$ be a braided category which is left rigid
(the definition of a left rigid category can be found in
\cite[XIV.2]{k} or \cite[IX.3]{maj}). If $V$ is an object
of ${\cal C}$ and $ev_V$ and $coev_V$ are the evaluation
and coevaluation maps associated to V, then following \cite{m2}
we define the representation-theoretic rank (or quantum dimension)
of $V$ as follows:
\[
\underline{\rm dim}(V)=ev _V\circ c _{V, V^*}\circ coev _V.
\]
${\;\;\;}$
If $H$ is a quasi-Hopf algebra then
the category ${}_H\Mm^{\rm fd}$ of finite dimensional modules over 
$H$ is left rigid. For $V\in {}_H\Mm$, its left dual is $V^*=\Hom(V,k)$, 
with left $H$-action $\lan h\cdot \varphi,v\ran= 
\lan\varphi,S(h)\cd v\ran$. The evaluation and coevaluation maps
are given for all $\v \in V^*$ and $v\in V$ by
\begin{eqnarray}
&&\hspace*{-2cm}
\ev_V(\varphi \ot v)=\varphi (\a \cd v),
\mbox{${\;\;\;}$}
\coev_V(1)=\sum \limits _i\b \cd v_i\ot v^i,\label{qrig}
\end{eqnarray}
where $\{v_i\}_i$ is a basis in $V$ with dual basis $\{v^i\}_i$ in $V^*$.

Therefore, if $H$ is a quasi-triangular quasi-Hopf algebra and
$V$ a finite dimensional left $H$-module it makes sense to consider
the representation-theoretic rank of $V$. If $R=R^1\ot R^2$ is an $R$-matrix
for $H$ then by \cite{bpv3} we have that
\begin{equation}\label{qd}
\un{\rm dim}(V)=\sum \limits _iv^i(S(R^2)\a R^1\b \cd v_i)=
{\rm Tr}(\eta),
\end{equation}
where $\eta :=S(R^2)\a R^1\b$. (Here ${\rm Tr}(\eta )$ is the trace of the
linear endomorphism of $V$ defined by $v\mapsto \eta \cd v$.)\\
${\;\;\;}$
Let $u$ be the element defined in (\ref{elmu}). By \cite{ac, bn3} we have
that $S(R^2)\a R^1=S(\a)u$, so by (\ref{sqina}) we obtain
\begin{equation}\label{eta}
\eta =S(S(\b)\a)u=u\smi (\a )\b .
\end{equation}
${\;\;\;}$
In the rest of this section $H$ will be a finite dimensional quasi-Hopf algebra,
and $\{e_i\}_{i=\ov{1, n}}$ a basis in $H$ with dual basis $\{e^i\}_{i=\ov{1, n}}$ in $H^*$. 
Our goal is to compute $\un{\rm dim}(H)$ and $\un{\rm dim}(D(H))$ within
the braided rigid category ${}_{D(H)}{\cal M}^{\rm fd}$. To this end
we shall compute for $D(H)$ the corresponding elements $u$ and $\eta$, denoted
in what follows by $u_D$ and $\eta _D$, respectively.

\begin{proposition}\label{pr4.1}
Let $H$ be a finite dimensional quasi-Hopf algebra, and $u_D$ and $\eta_D$ the
corresponding elements $u$ and $\eta$ for $D(H)$, the quantum double of $H$. Then
\begin{equation}\label{ue}
u_D=\sum \limits _{i=1}^n\b \rh \ov{S}^{-1}(e^i)\bowtie e_i~~\mbox{and}~~
\eta _D=\sum \limits _{i=1}^n\b \rh \ov{S}^{-1}(e^i)\bowtie e_i\smi (\a)\b .
\end{equation}
\end{proposition}
\begin{proof}
Let us start by noting that (\ref{pf}), (\ref{fgab}) and (\ref{q5}) imply
\begin{equation}\label{f5}
f^1_1p^1\ot f^1_2p^2S(f^2)=g^1S(\tqlb)\ot g^2S(\tqla).
\end{equation}
Secondly, observe that the definition (\ref{anddf}) of the antipode $S_D$ of
$D(H)$ can be reformulated as follows:
\begin{eqnarray*}
&&\hspace*{-1.5cm}
S_D(\v \bowtie h)\\
&{{\rm (\ref{mdd})}\atop =}&
\left(\va \bowtie S(h)\right)\left((f^1_1p^1)_1U^1\rh \ov{S}^{-1}(\v)\lh
f^2\smi (f^1_2p^2)\bowtie (f^1_1p^1)_2U^2\right)\\
&{{\rm (\ref{f5})}\atop =}&
\left(\va \bowtie S(h)\right)
\left(g^1_1S(\tqlb)_1U^1\rh \ov{S}^{-1}(\v)\lh \tqla \smi (g^2)
\bowtie g^1_2S(\tqlb)_2U^2\right)\\
&{{\rm (\ref{U}, \ref{ca})}\atop =}&
\left(\va \bowtie S(h)\right)\left(g^1_1G^1S(q^2\tilde{q}^2_2)\rh \ov{S}^{-1}(\v )\lh
\tqla \smi (g^2)\bowtie g^1_2G^2S(q^1\tilde{q}^2_1)\right),
\end{eqnarray*}
where we denoted by $G^1\ot G^2$ another copy of $f^{-1}$. Now, we claim that
\begin{equation}\label{f6}
S_D(\cal{R}^2)\a _D{\cal R}^1=\sum \limits_{i=1}^n
\b \rh \ov{S}^{-1}(e^i)\lh \a \bowtie e_i,
\end{equation}
where $R_D={\cal R}^1\ot {\cal R}^2$ is the $R$-matrix of
$D(H)$ defined in (\ref{rdd}). Indeed, we can easily check that
\begin{equation}\label{f7}
\ov{S}^{-1}(h\rh \v)=\ov{S}^{-1}(\v)\lh S(h)~~{\rm and}~~
\ov{S}^{-1}(\v \lh h)=S(h)\rh \ov{S}^{-1}(\v),
\end{equation}
for all $\v \in H^*$ and $h\in H$. Now, we calculate:
\begin{eqnarray*}
&&\hspace*{-1.5cm}
S_D(\cal{R}^2)\a _D{\cal R}^1\\
&{{\rm (\ref{rdd}, \ref{abdd})}\atop =}&\sum \limits_{i=1}^n
S_D(e^i\bowtie p^1_2)(\va \bowtie \a )(\va \bowtie \smi (p^2)e_ip^1_1)\\
&{{\rm (\ref{mdd})}\atop =}&\sum \limits_{i=1}^n
\left(S(p^1_2)_1g^1\right)_1G^1S(q^2\tilde{q}^2_2)\rh \ov{S}^{-1}(e^i)\lh
\tqla\smi(S(p^1_2)_2g^2)\\
&&\hspace*{1cm}
\bowtie \left(S(p^1_2)_1g^1\right)_2G^2S(q^1\tilde{q}^2_1)\a\smi (p^2)e_ip^1_1\\
&{{\rm (\ref{ca})}\atop =}&\sum \limits_{i=1}^n
g^1_1G^1S\left(q^2(\tqlb p^1_{(2, 2)})_2\right)\rh \ov{S}^{-1}(p^1_1\rh e^i)
\lh \tqla p^1_{(2, 1)}\smi (g^2)\\
&&\hspace*{1cm}
\bowtie g^1_2G^2S\left(q^1(\tqlb p^1_{(2, 2)})_1\right)\a \smi (p^2)e_i\\
&{{\rm (\ref{f7}, \ref{ql1a})}\atop =}&\sum \limits_{i=1}^n
g^1_1G^1S(q^2p^1_2\tilde{q}^2_2)\rh \ov{S}^{-1}(e^i\lh \smi (p^2))\lh
\tqla \smi (g^2)\\
&&\hspace*{1cm}
\bowtie g^1_2G^2S(q^1p^1_1\tilde{q}^2_1)\a e_i\\
&{{\rm (\ref{f7}, \ref{pqra})}\atop =}&\sum \limits_{i=1}^n
g^1_1G^1S(\tilde{q}^2_2)\rh \ov{S}^{-1}(e^i)\lh \tqla \smi (g^2)
\bowtie g^1_2G^2S(\tilde{q}^2_1)\a e_i\\
&{{\rm (\ref{f7})}\atop =}&\sum \limits_{i=1}^n
g^1_1G^1\rh \ov{S}^{-1}(e^i)\lh \tqla \smi (g^2)
\bowtie g^1_2G^2S(\tilde{q}^2_1)\a \tilde{q}^2_2e_i\\
&{{\rm (\ref{q5}, \ref{ql})}\atop =}&\sum \limits_{i=1}^n
e^i\bowtie g^1_2G^2\a \smi (\a \smi (g^2)e_ig^1_1G^1)\\
&{{\rm (\ref{fgab}, \ref{q5})}\atop =}&\sum \limits_{i=1}^n
e^i\bowtie \smi (\a e_i\b )=\sum \limits_{i=1}^n
\b \rh \ov{S}^{-1}(e^i)\lh \a \bowtie e_i.
\end{eqnarray*}
We are now able to calculate the element $u_D$. Since $H$ can
be viewed as a quasi-Hopf subalgebra of $D(H)$ via the morphism $i_D$ it
follows that the corresponding element $p_R$ for $D(H)$ is
$(p_R)_D=p^1_D\ot p^2_D=\va \bowtie p^1\ot \va \bowtie p^2$. Therefore:
\begin{eqnarray*}
u_D&{{\rm (\ref{elmu})}\atop =}&
S_D({\cal R}^2p^2_D)\a _D{\cal R}^1p^1_D=(\va \bowtie S(p^2))S_D({\cal R}^2)\a _D
{\cal R}^1(\va \bowtie p^1)\\
&{{\rm (\ref{f6}, \ref{mdd})}\atop =}&\sum \limits _{i=1}^n
S(p^2)_{(1, 1)}\b \rh \ov{S}^{-1}(e^i)\lh \a \smi (S(p^2)_2)
\bowtie S(p^2)_{(1, 2)}e_ip^1\\
&=&\sum \limits _{i=1}^n
e^i\bowtie S(p^2)_{(1, 2)}\smi \left(\a \smi (S(p^2)_2)e_iS(p^2)_{(1, 1)}\b \right)p^1\\
&{{\rm (\ref{q5})}\atop =}&\sum \limits _{i=1}^n
e^i\bowtie \smi (S(p^1)\a p^2e_i\b )\\
&{{\rm (\ref{qr}, \ref{q6})}\atop =}&\sum \limits _{i=1}^n
e^i\bowtie \smi (e_i\b)=\sum \limits _{i=1}^n\b \rh \ov{S}^{-1}(e^i)\bowtie e_i,
\end{eqnarray*}
as claimed. It is clear now that the above equality and (\ref{eta}) imply
the expression of $\eta _D$ in (\ref{ue}), so our proof is complete.
\end{proof}

\subsection{The representation-theoretic rank of $H$}
We start to compute the repre-sentation-theoretic rank (or
quantum dimension) of $H$ within the braided rigid category
${}_{D(H)}{\cal M}^{\rm fd}$. Let us start by noting that the action
$\act$ obtained in (\ref{aadd}) can be rewritten as follows:
\begin{eqnarray*}
&&\hspace*{-2cm}
(\v \bowtie h)\act h'\\
&{{\rm (\ref{q3}, \ref{q5})}\atop =}&
\< \v , \smi (Y^3)q^2Y^2_2y^3_2\smi (\tqla y^2_1(h\tr h')_1g^1)y^1\>\\
&&\hspace*{3cm}
Y^1\tqlb y^2_2(h\tr h')_2g^2S(q^1Y^2_1y^3_1)\\
&{{\rm (\ref{ca}, \ref{U})}\atop =}&
\< \v , \smi \left(\tqla (y^2(h\tr h')S(Y^2y^3))_1U^1Y^3\right)y^1\>\\
&&\hspace*{3cm}
Y^1\tqlb (y^2(h\tr h')S(Y^2y^3))_2U^2\\
&{{\rm (\ref{ql1a})}\atop =}&
\< \v , \smi \left(\tqla (Y^1_2y^2(h\tr h')S(Y^2y^3))_1U^1Y^3\right)Y^1_1y^1\>\\
&&\hspace*{3cm}
\tqlb (Y^1_2y^2(h\tr h')S(Y^2y^3))_2U^2.
\end{eqnarray*}
Hence we have showed that for all $\v \in H^*$ and $h, h'\in H$ we have
\begin{eqnarray}
&&\hspace*{-2cm}
(\v \bowtie h)\act h'=
\< \v , \smi \left(\tqla (Y^1_2y^2(h\tr h')S(Y^2y^3))_1U^1Y^3\right)
Y^1_1y^1\>\nonumber\\
&&\hspace*{4cm}
\tqlb (Y^1_2y^2(h\tr h')S(Y^2y^3))_2U^2.\label{act2}
\end{eqnarray}
So this action defines on $H$ a left $D(H)$-module structure, and on
$H_0$ a left $D(H)$-module algebra structure.\\
${\;\;\;}$
In order to "simplify" the computation for $\un{\rm dim}(H)$
we need the following formulas.
\begin{lemma}\label{lm4.2}
Let $H$ be a finite dimensional quasi-Hopf algebra and
$\{e_i\}_i$ a basis in $H$ with
dual basis $\{e^i\}$. Then for all $h, h', h{''}\in H$
the following relations hold:
\begin{eqnarray}
&&\sum \limits _{i=1}^n\< e^i, \smi (\b )S^{-2}(\tQla (e_i)_1h')h\tqlb
\tilde{Q}^2_2(e_i)_{(2, 2)}h{''}
\smi (\tqla \tilde{Q}^2_1(e_i)_{(2, 1)})\>\nonumber\\
&&\hspace*{1cm}
=\sum \limits _{i=1}^n\< e^i, \smi (\b )S^{-2}(\tQla (e_i)_1h')
\tqlb \tilde{Q}^2_2(e_i)_{(2, 2)}h{''}
\smi (\tqla \tilde{Q}^2_1(e_i)_{(2, 1)})h\>,\label{f8}
\end{eqnarray}
\begin{eqnarray}
&&\sum \limits _{i=1}^n\< e^i, \smi (\b )S^{-2}(\tQla (e_i)_1X^1p^1_1h')h_1\tqlb
\tilde{Q}^2_2(e_i)_{(2, 2)}X^3p^2S(h_2)h{''}\nonumber\\
&&\hspace*{1cm}\times
\smi (\tqla \tilde{Q}^2_1(e_i)_{(2, 1)}X^2p^1_2)\>=\sum \limits _{i=1}^n
< e^i, \smi (\b )S^{-2}(\tQla (e_i)_1X^1p^1_1h_1h')\nonumber\\
&&\hspace*{1cm}\times
\tqlb \tilde{Q}^2_2(e_i)_{(2, 2)}X^3p^2h{''}
\smi (\tqla \tilde{Q}^2_1(e_i)_{(2, 1)}X^2p^1_2h_2)\>,\label{f9}
\end{eqnarray}
where we denoted $q_L=\tqla \ot \tqlb =\tQla \ot \tQlb$ and $p_R=p^1\ot p^2$.
\end{lemma}
\begin{proof}
In order to prove (\ref{f8}) we shall apply (\ref{ql1a}) twice, and then
the properties of dual bases and (\ref{q5}). Explicitly,
\begin{eqnarray*}
&&\hspace*{-1.5cm}
\sum \limits _{i=1}^n\< e^i, \smi (\b )S^{-2}(\tQla (e_i)_1h')h\tqlb
\tilde{Q}^2_2(e_i)_{(2, 2)}h{''}
\smi (\tqla \tilde{Q}^2_1(e_i)_{(2, 1)})\>\\
&=&\sum \limits _{i=1}^n\< e^i, \smi (\b )S^{-2}(\tQla (e_i)_1h')\tqlb (h_2
\tilde{Q}^2)_2(e_i)_{(2, 2)}h{''}\\
&&\hspace*{2cm}\times
\smi (\tqla (h_2\tilde{Q}^2)_1(e_i)_{(2, 1)})h_1\>
\end{eqnarray*}
\begin{eqnarray*}
&=&\sum \limits _{i=1}^n\< e^i, \smi (\b )
S^{-2}(S(h_{(2, 1)})\tQla (h_{(2, 2)}e_i)_1h')\tqlb
\tilde{Q}^2_2(h_{(2, 2)}e_i)_{(2, 2)}h{''}\\
&&\hspace*{2cm}\times
\smi (\tqla \tilde{Q}^2_1(h_{(2, 2)}e_i)_{(2, 1)})h_1\>\\
&=&\sum \limits _{i=1}^n\< e^i, h_{(2, 2)}
\smi (h_{(2, 1)}\b )S^{-2}(\tQla (e_i)_1h')\tqlb
\tilde{Q}^2_2(e_i)_{(2, 2)}h{''}\\
&&\hspace*{2cm}\times
\smi (\tqla \tilde{Q}^2_1(e_i)_{(2, 1)})h_1\>\\
&=&\sum \limits _{i=1}^n\< e^i,
\smi (\b )S^{-2}(\tQla (e_i)_1h')
\tqlb \tilde{Q}^2_2(e_i)_{(2, 2)}h{''}
\smi (\tqla \tilde{Q}^2_1(e_i)_{(2, 1)})h\>.
\end{eqnarray*}
In a similar manner we can prove (\ref{f9}). It follows applying
(\ref{f8}), dual basis, (\ref{q1}) and (\ref{qr1}), we leave the
details to the reader.
\end{proof}
Now, equation (\ref{gd}) shows by using (\ref{q3}) and (\ref{q5}) that
\begin{eqnarray}
&&\gamma =\gamma ^1\ot \gamma ^2=S(x^1X^2)\a x^2X^3_1\ot S(X^1)\a x^3X^3_2,\label{g2}\\
&&\delta =\delta ^1\ot \delta ^2=x^1\b S(x^3_2X^3)\ot x^2X^1\b S(x^3_1X^2).\label{d2}
\end{eqnarray}
We finally need the following formula
\begin{equation}\label{U2}
p_R=\Delta (S(\tpla))U(\tplb \ot 1),
\end{equation}
which can be found in \cite{hn3}. We are now able to compute $\un{\rm dim}(H)$.

\begin{proposition}\label{pr4.3}
Let $H$ be a finite dimensional quasi-Hopf algebra. Then the representation-theoretic
rank of $H$ is
\begin{equation}\label{qdH}
\un{\rm dim}(H)={\rm Tr}\left(h\mapsto S^{-2}(S(\b )\a h\b S(\a))\right).
\end{equation}
\end{proposition}
\begin{proof}
We know from Proposition \ref{pr4.1}
that in the quantum double case the element $\eta _D$ is given by
\[
\eta _D=\sum \limits _{i=1}^n\b \rh \ov{S}^{-1}(e^i)\bowtie e_i\smi (\a)\b
=\sum \limits _{i=1}^ne^i\bowtie \smi (\a e_i\b )\b .
\]
We set $p_R=p^1\ot p^2=P^1\ot P^2$, $q_L=\tqla \ot \tqlb=\tQla \ot \tQlb$ and
$f=f^1\ot f^2=F^1\ot F^2$.
Then by (\ref{qd}) and the above expression of $\eta _D$ we have:
\begin{eqnarray*}
&&\hspace*{-2cm}\un{\rm dim}(H)\\
&=&\sum \limits _{i, j=1}^n
\< e^j, \left(e^i\bowtie \smi (\a e_i\b)\b\right)\act e_j\>\\
&{{\rm (\ref{act2})}\atop =}&\sum \limits _{i, j=1}^n
\< e^i, \smi \left(\tqla (Y^1_2y^2(\smi (\a e_i\b)\b \tr e_j)
S(Y^2y^3))_1U^1Y^3\right)Y^1_1y^1\>\\
&&\hspace*{5mm}
\< e^j, \tqlb (Y^1_2y^2(\smi (\a e_i\b)\b \tr e_j)S(Y^2y^3))_2U^2\>\\
&=&\sum \limits _{i, j, k=1}^n
\<e^k, Y^1_2y^2(\smi (\a e_i\b)\b \tr e_j)S(Y^2y^3)\>\\
&&\hspace*{5mm}
\<Y^1_1y^1\rh e^i, \smi (\tqla (e_k)_1U^1Y^3)\>
\<e^j , \tqlb (e_k)_2U^2\>
\end{eqnarray*}
\begin{eqnarray*}
&=&\sum \limits _{i, k=1}^n
\< e^k , Y^1_2y^2\left(\smi (\a e_iY^1_1y^1\b )
\b \tr \tqlb (e_k)_2U^2\right)S(Y^2y^3)\>\\
&&\hspace*{5mm}
\< e^i, \smi (\tqla (e_k)_1U^1Y^3)\>\\
&{{\rm (\ref{lya2}, \ref{ca}, \ref{gdf})}\atop =}&\sum \limits _{i, k=1}^n
\< e^k, Y^1_2y^2\smi (f^2Y^1_{(1, 2)}y^1_2\delta ^2)
\left(\smi (\a e_i)\b \tr \tqlb (e_k)_2U^2\right)\\
&&\hspace*{5mm}
\times f^1Y^1_{(1, 1)}y^1_1\delta ^1S(Y^2y^3)\>
\< e^i, \smi (\tqla (e_k)_1U^1Y^3)\>\\
&{{\rm (\ref{d2}, \ref{q3}, \ref{qr})}\atop =}&\sum \limits _{i, k=1}^n
\< e^k, Y^1_2\smi (f^2Y^1_{(1, 2)}p^2)\left(\smi (\a e_i)\b \tr \tqlb (e_k)_2U^2\right)\\
&&\hspace*{5mm}\times
f^1Y^1_{(1, 1)}p^1\b S(Y^2)\> \< e^i, \smi (\tqla (e_k)_1U^1Y^3)\>\\
&{{\rm (\ref{qr1}, \ref{ql})}\atop =}&\sum \limits _{i, k=1}^n
\< S(\tpla)\rh e^k, \smi (f^2p^2)\left(\smi (\a e_i)\b \tr
\tqlb (e_k)_2U^2\right)f^1p^1\>\\
&&\hspace*{5mm}
\< e^i, \smi (\tqla (e_k)_1U^1\tplb )\>\\
&=&\sum \limits _{i, k=1}^n
\< e^k, \smi (f^2p^2)\left(\smi (\a e_i)\b \tr
\tqlb (e_k)_2S(\tpla)_2U^2\right)f^1p^1\>\\
&&\hspace*{5mm}
\< e^i, \smi (\tqla (e_k)_1S(\tpla)_1U^1\tplb )\>\\
&{{\rm (\ref{U2})}\atop =}&\sum \limits _{k=1}^n
\< e^k, \smi (f^2p^2)\left(\smi (\a \smi (\tqla (e_k)_1P^1))\b \tr
\tqlb (e_k)_2P^2\right)f^1p^1\>\\
&{{\rm (\ref{lya2}, \ref{ca}, \ref{qr})}\atop =}&\sum \limits _{k=1}^n
\< e^k\lh x^3, \smi (f^2\smi (F^1\tilde{q}^1_1(e_k)_{(1, 1)}p^1_1g^1)x^2\b )\\
&&\hspace*{5mm}\times
\left(\smi (\a)\b \tr \tqlb (e_k)_2p^2\right)f^1
\smi (F^2\tilde{q}^1_2(e_k)_{(1, 2)}p^1_2g^2)x^1\>\\
&=&\sum \limits _{k=1}^n
\< e^k, \smi (f^2\smi (F^1\tilde{q}^1_1x^3_{(1, 1)}(e_k)_{(1, 1)}p^1_1g^1)x^2\b )\\
&&\hspace*{5mm}\times
\left(\smi (\a)\b \tr \tqlb x^3_2(e_k)_2p^2\right)f^1
\smi (F^2\tilde{q}^1_2x^3_{(1, 2)}(e_k)_{(1, 2)}p^1_2g^2)x^1\>\\
&{{\rm (\ref{ql2}, \ref{ca}, \ref{gdf})}\atop =}&\sum \limits _{k=1}^n
\< e^k, \smi (\gamma ^2\smi (\tQla X^1(e_k)_{(1, 1)}p^1_1g^1)\b )\\
&&\hspace*{5mm}\times
\b _1\tqlb \tilde{Q}^2_2X^3(e_k)_2p^2S(\b _2)\gamma ^1
\smi (\tqla \tilde{Q}^2_1X^2(e_k)_{(1, 2)}p^1_2g^2)\>\\
&{{\rm (\ref{q1}, \ref{f8})}\atop =}&\sum \limits _{k=1}^n
\< e^k, \smi (\gamma ^2\smi (\tQla (e_k)_1X^1p^1_1g^1)\b )\\
&&\hspace*{5mm}\times
\tqlb \tilde{Q}^2_2(e_k)_{(2, 2)}X^3p^2S(\b _2)\gamma ^1
\smi (\tqla \tilde{Q}^2_1(e_k)_{(2, 1)}X^2p^1_2g^2)\b _1\>\\
&{{\rm (\ref{f9}, \ref{gdf})}\atop =}&\sum \limits _{k=1}^n
\< e^k, \smi (\gamma ^2\smi (\tQla (e_k)_1X^1p^1_1\delta ^1)\b )\\
&&\hspace*{5mm}\times
\tqlb \tilde{Q}^2_2(e_k)_{(2, 2)}X^3p^2\gamma ^1
\smi (\tqla \tilde{Q}^2_1(e_k)_{(2, 1)}X^2p^1_2\delta ^2)\>\\
&{{\rm (\ref{gd}, \ref{g2})}\atop =}&\sum \limits _{k=1}^n
\< e^k, \smi (\b )S^{-2}(\tQla (e_k)_1X^1p^1_1Y^1_1x^1\b S(S(Z^1)\a y^3Z^3_2Y^3))\\
&&\hspace*{5mm}\times
\tqlb \tilde{Q}^2_2(e_k)_{(2, 2)}X^3p^2\smi (\tqla \tilde{Q}^2_1(e_k)_{(2, 1)}
X^2p^1_2Y^1_2
\end{eqnarray*}
\begin{eqnarray*}
&&\hspace*{5mm}\times
x^2\b S(S(y^1Z^2)\a y^2Z^3_1Y^2x^3)\>\\
&{{\rm (\ref{qr1}, \ref{q1})}\atop =}&\sum \limits _{k=1}^n
\< Y^1_1\rh e^k, \smi (\b )S^{-2}(\tQla (e_k)_1X^1p^1_1x^1\b S(S(Z^1)\a y^3Z^3_2Y^3))\\
&&\hspace*{5mm}\times
\tqlb \tilde{Q}^2_2(e_k)_{(2, 2)}X^3p^2\smi (\tqla \tilde{Q}^2_1(e_k)_{(2, 1)}
X^2p^1_2\\
&&\hspace*{5mm}\times
x^2\b S(S(y^1Z^2Y^1_2)\a y^2Z^3_1Y^2x^3)\>\\
&{{\rm (\ref{f8}, \ref{q3}, \ref{q5})}\atop =}&\sum \limits _{k=1}^n
\< e^k, \smi (\b )S^{-2}(\tQla (e_k)_1X^1p^1_1x^1\b S(\a))Y^1
\tqlb \tilde{Q}^2_2(e_k)_{(2, 2)}X^3p^2\\
&&\hspace*{5mm}\times
\smi (\tqla \tilde{Q}^2_1(e_k)_{(2, 1)}X^2p^1_2x^2\b S(S(Y^2)\a Y^3x^3))\>\\
&{{\rm (\ref{qr}, \ref{f8})}\atop =}&\sum \limits _{k=1}^n
\< q^1\rh e^k, \smi (\b)S^{-2}(\tQla (e_k)_1X^1p^1_1P^1\b S(\a))\tqlb \tilde{Q}^2_2\\
&&\hspace*{5mm}\times
(e_k)_{(2, 2)}X^3p^2S(q^2)\smi (\tqla \tilde{Q}^2_1(e_k)_{(2, 1)}X^2p^1_2P^2)\>\\
&{{\rm (\ref{q1})}\atop =}&\sum \limits _{k=1}^n
\< e^k, \smi (\b)S^{-2}(\tQla (e_k)_1X^1(q^1_1p^1)_1P^1\b S(\a))\tqlb \tilde{Q}^2_2\\
&&\hspace*{5mm}\times
(e_k)_{(2, 2)}X^3q^1_2p^2S(q^2)\smi (\tqla \tilde{Q}^2_1(e_k)_{(2, 1)}
X^2(q^1_1p^1)_2P^2)\>\\
&{{\rm (\ref{pqr}, \ref{qr})}\atop =}&\sum \limits _{k=1}^n
\< e^k, \smi (\b)S^{-2}(\tQla (e_k)_1\b S(\a))\tqlb \\
&&\hspace*{5mm}\times
\tilde{Q}^2_2(e_k)_{(2, 2)}\smi (\tqla \tilde{Q}^2_1(e_k)_{(2, 1)}\b )\>\\
&{{\rm (\ref{q5}, \ref{ql})}\atop =}&\sum \limits _{k=1}^n
\< e^k, \smi (\b)S^{-2}(\a e_k\b S(\a))\tqlb \smi (\tqla \b )\>\\
&{{\rm (\ref{ql}, \ref{q6})}\atop =}&\sum \limits _{k=1}^n
\< e^k, S^{-2}(S(\b)\a e_k\b S(\a ))\>
={\rm Tr}\left(h\mapsto S^{-2}(S(\b)\a h\b S(\a ))\right),
\end{eqnarray*}
so the proof is finished.
\end{proof}
${\;\;\;}$
In \seref{5} we will see that the representation-theoretic rank of $H$
can be expressed in terms of integrals in $H$ and $H^*$. This result is
strictly connected to the trace formula for quasi-Hopf algebras.
\subsection{The representation-theoretic rank of $D(H)$}
Let $H$ be a finite dimensional quasi-triangular quasi-Hopf algebra. Then
$H$ is an object in its own category of finite dimensional
representations via the left regular action, so it makes sense to
consider $\un{\rm dim}(H)$.\\
${\;\;\;}$
The purpose of this subsection is to compute $\un{\rm dim}(D(H))$.
As we will see the computation is harder than the one for
$\un{\rm dim}(H)$ but the result will be the same. Again, we need some
preliminary work.\\
${\;\;\;}$
Recall that $t\in H$ is called a left (respectively right) integral
in $H$ if $ht=\va(h)t$ (respectively $th=\va(h)t$), for all $h\in H$. We denote
by $\int _l$ and $\int _r$ the space of left and right integrals in $H$.
When $H$ is finite dimensional we have that
${\rm dim}(\int _l)={\rm dim}(\int _r)=1$, $S(\int _l)=\int _r$
and $S(\int _r)=\int _l$ (see \cite{hn3, bc1}). In addition, if we define
\[
\mf{P}(h)=\sum \limits _{i=1}^n \< e^i, \b S^2(q^2(e_i)_2)h\>q^1(e_i)_1,
~~\forall~~h\in H,
\]
then by \cite{pv2} we have that $\mf{P}(h)\in \int _l$, for all $h\in H$, and
$\mf{P}(t)=t$ for any $t\in \int _l$. Therefore, $\mf{P}$ defines a projection from
$H$ to $\int _l$. Replacing the quasi-Hopf algebra $H$ by $H^{cop}$ we
obtain a second projection onto the space of
left integrals, denoted in what follows by $\tilde{\mf{P}}$. Since in $H^{cop}$
we have $(q_R)_{\rm cop}=\tqlb \ot \tqla$ we obtain
\begin{equation}\label{f10}
\tilde{\mf{P}}(h)=\sum \limits _{i=1}^n\<e^i, \smi (\b)S^{-2}(\tqla (e_i)_1)h\>
\tqlb (e_i)_2\in \int _l, ~~\forall~~h\in H.
\end{equation}
${\;\;\;}$
We finally need the following formulas.
\begin{lemma}\label{lm4.4}
In a quasi-Hopf algebra $H$ the following relations hold:
\begin{eqnarray}
&&\Omega ^1_1\d ^1S^2(\Omega ^4)\ot \Omega ^1_{(2, 1)}\delta ^2_1g^1S(\Omega ^3)
\ot \Omega ^1_{(2, 2)}\delta ^2_2g^2S(\Omega ^2)\ot \Omega ^5\nonumber\\
&&\hspace*{2cm}
=X^1p^1_1P^1S(f^1\tpla )\ot X^2p^1_2P^2\ot X^3p^2
\ot \smi (f^2\tplb),\label{f11}\\
&&\gamma ^1X^1\ot f^1\gamma ^2_1X^2\ot f^2
\gamma ^2_2X^3=S(X^3)f^1\gamma ^1_1\ot
S(X^2)f^2\gamma ^1_2\ot S(X^1)\gamma ^2,\label{f12}\\
&&q^1x^1_1\ot \smi (x^2)q^2x^1_2\ot x^3=
X^1\ot \smi (\tqla X^3_1)X^2\ot \tqlb X^3_2.\label{f13}
\end{eqnarray}
Here $\Omega =\Omega ^1\ot \cdots \ot \Omega ^5$, $\d =\d ^1\ot \d ^2$,
$\gamma =\gamma ^1\ot \gamma ^2$, $f=f^1\ot f^2$, $f^{-1}=g^1\ot g^2$,
$q_R=q^1\ot q^2$, $p_R=p^1\ot p^2=P^1\ot P^2$
and $q_L=\tqla \ot \tqlb$ are the elements defined in
(\ref{O}), (\ref{gd}), (\ref{f}), (\ref{g}), (\ref{qr}) and
(\ref{ql}),respectively.
\end{lemma}
\begin{proof}
Observe that the element $\delta $ in (\ref{gd}) can be rewritten as
\[
\d =Y^1_1p^1\b S(Y^3)\ot Y^1_2p^2S(Y^2).
\]
Now, using the above description for $\d$ and the definition of $\Omega$
we compute:
\begin{eqnarray*}
&&\hspace*{-2.2cm}
\Omega ^1_1\d ^1S^2(\Omega ^4)\ot \Omega ^1_{(2, 1)}\delta ^2_1g^1S(\Omega ^3)
\ot \Omega ^1_{(2, 2)}\delta ^2_2g^2S(\Omega ^2)\ot \Omega ^5\\
&{{\rm (\ref{ca})}\atop =}&
X^1_{(1, 1)_1}y^1_1p^1\b S(f^1X^2)\ot X^1_{(1, 1)_{(2, 1)}}y^1_{(2, 1)}p^2_1
g^1S(X^1_2y^3)\\
&&\hspace*{1cm}
\ot X^1_{(1, 1)_{(2, 2)}}y^1_{(2, 2)}p^2_2
g^2S(X^1_{(1, 2)}y^2)\ot \smi (f^2X^3)\\
&{{\rm (\ref{pr}, \ref{q1})}\atop =}&
Y^1\left((X^1_1)_{(1, 1)}p^1\right)_1P^1\b S(f^1X^2)\ot Y^2
\left((X^1_1)_{(1, 1)}p^1\right)_2P^2S(X^1_2)\\
&&\hspace*{1cm}
\ot Y^3(X^1_1)_{(1, 2)}p^2S((X^1_1)_2)\ot \smi (f^2X^3)\\
&{{\rm (\ref{qr1}, \ref{ql})}\atop =}&
Y^1p^1_1P^1S(f^1\tpla)\ot Y^2p^1_2P^2\ot Y^3p^2\ot \smi (f^2\tplb),
\end{eqnarray*}
so the equality in (\ref{f11}) is proved. The relation in (\ref{f12})
follows more easily since
\begin{eqnarray*}
&&\hspace*{-2cm}
\gamma ^1X^1\ot f^1\gamma ^2_1X^2\ot f^2
\gamma ^2_2X^3\\
&{{\rm (\ref{gdf})}\atop =}&
F^1\a _1X^1\ot f^1F^2_1\a _{(2, 1)}X^2\ot f^2F^2_2\a _{(2, 2)}X^3\\
&{{\rm (\ref{q1}, \ref{pf})}\atop =}&
S(X^3)f^1F^1_1\a _{(1, 1)}\ot S(X^2)f^2F^1_2\a _{(1, 2)}\ot S(X^1)F^2\a _2\\
&{{\rm (\ref{gdf})}\atop =}&
S(X^3)f^1\gamma ^1_1\ot S(X^2)f^2\gamma ^1_2\ot S(X^1)\gamma ^2,
\end{eqnarray*}
where we denoted by $F^1\ot F^2$ another copy of $f$. Finally, (\ref{f13}) is an
immediate consequence of (\ref{q3}) and (\ref{q5}).
\end{proof}
${\;\;\;}$
We can now compute the representation-theoretic rank of $D(H)$. The next
result generalizes \cite[Proposition 2.1]{m2}.
\begin{proposition}\label{pr4.5}
Let $H$ be a finite dimensional quasi-Hopf algebra and $D(H)$ its
quantum double. Then
\[
\un{\rm dim}(D(H))=\un{\rm dim }(H)={\rm Tr}\left(h\mapsto S^{-2}
(S(\b)\a h\b S(\a))\right).
\]
\end{proposition}
\begin{proof}
We set $p_R=p^1\ot p^2=P^1\ot P^2$,
$q_R=q^1\ot q^2=Q^1\ot Q^2$ and
$f=f^1\ot f^2=F^1\ot F^2={\cal F}^1\ot {\cal F}^2$. In what follows,
we shall not perform all the computations but we shall point out
the relations which are used in every step.\\
${\;\;\;}$
The expression of $\eta_D$ in Proposition \ref{pr4.1} allows us to compute:
${\;\;}$\vspace*{-1cm}
\begin{eqnarray*}
&&\hspace*{-2cm}
\un{\rm dim}(D(H))=\sum \limits _{i, j=1}^n\< e_i\bowtie e^j, \eta _D(e^i\bowtie e_j)\>\\
&=&\sum \limits _{i, j, k=1}^n\<e_i\bowtie e^j,
(\b \rh \ov{S}^{-1}(e^k)\lh S(\b)\a \bowtie e_k)
(e^i\bowtie e_j)\>\\
&{{\rm (\ref{mdd})}\atop =}&\sum \limits _{i, j, k=1}^n
\<\ov{S}^{-1}(e^k), S(\b)\a \Omega ^5(e_i)_1\Omega ^1\b \>\\
&&\hspace*{1cm}
\< e^i, \smi ((e_k)_2)\Omega ^4(e_i)_2\Omega ^2(e_k)_{(1, 1)}\>
\< e^j, \Omega ^3(e_k)_{(1, 2)}e_j\>\\
&{{\rm (\ref{ca}, \ref{gdf})}\atop =}&\sum \limits _{i, j=1}^n
\< e^j, \Omega ^3
\smi \left(S(\b _1)\gamma ^2\Omega ^5_2(e_i)_{(1, 2)}
\Omega ^1_2\d ^2\right)_2e_j\>\\
&&\hspace*{1cm}
\< e^i, S^{-2}(S(\b _2)\gamma ^1\Omega ^5_1(e_i)_{(1, 1)}\Omega ^1_1\d ^1)
\Omega ^4(e_i)_2\Omega ^2\\
&&\hspace*{1cm}\times
\smi \left(S(\b _1)\gamma ^2\Omega ^5_2(e_i)_{(1, 2)}
\Omega ^1_2\d ^2\right)_1\>\\
&{{\rm (\ref{ca}, \ref{f11})}\atop =}&\sum \limits _{i, j=1}^n
\< e^i, S^{-2}\left(S(\b _2)\gamma ^1\smi (F^2\tplb )_1
((e_i)_1)_ 1Y^1p^1_1P^1S(F^1\tpla )\right)(e_i)_2\\
&&\hspace*{1cm}\times
\smi \left(f^2S(\b _1)_2\gamma ^2_2
\smi (F^2\tplb )_{(2, 2)}((e_i)_1)_{(2 , 2)}Y^3p^2\right)\>\\
&&\hspace*{1cm}
\< e^j, \smi \left(f^1S(\b _1)_1\gamma ^2_1
\smi (F^2\tplb )_{(2, 1)}((e_i)_1)_{(2, 1)}Y^2p^1_2P^2\right)e_j\>\\
&{{\rm (\ref{q1}, \ref{qr1})}\atop =}&\sum \limits _{i, j=1}^n
\< e^i, S^{-2}\left(S(\b _2)\gamma ^1Y^1
(\smi (F^2\tplb)_1p^1)_1(e_i)_1P^1S(F^1\tpla )\right)\\
&&\hspace*{1cm}\times
\smi \left(f^2S(\b _1)_2\gamma ^2_2Y^3\smi (F^2\tplb )_2p^2\right)\>\\
&&\hspace*{1cm}
\< e^j, \smi \left(f^1S(\b _1)_1\gamma ^2_1Y^2
(\smi (F^2\tplb )_1p^1)_2(e_i)_2P^2\right)e_j\>\\
&{{\rm (\ref{ca}, \ref{qr}, \ref{pf})}\atop =}&\sum \limits _{i, j=1}^n
\< e^i, S^{-2}\left(S(\b _2)\gamma ^1Y^1\smi (F^2x^3\tilde{p}^2_2g^2)_1(e_i)_1
P^1S({\cal F}^1F^1_1x^1\tpla )\right)\\
&&\hspace*{1cm}\times
\smi \left(f^2S(\b _1)_2\gamma ^2_2Y^3
\smi ({\cal F}^2F^1_2x^2\tilde{p}^2_1g^1)\b \right)\>\\
&&\hspace*{1cm}
\< e^j, \smi \left(f^1S(\b _1)_1\gamma ^2_1Y^2\smi (F^2x^3\tilde{p}^2_2g^2)_2
(e_i)_2P^2\right)e_j\>\\
&{{\rm (\ref{fgab}, \ref{q5}, \ref{ql})}\atop =}&\sum \limits _{i, j=1}^n
\< e^i, S^{-2}(S(\b _2)\gamma ^1Y^1\smi (\tqlb \tilde{p}^2_2g^2)_1(e_i)_1P^1
S(\tpla)\tqla \tilde{p}^2_1g^1\\
&&\hspace*{1cm}\times
S(f^2S(\b _1)_2\gamma ^2_2Y^3))\>
\< e^j, \smi (f^1S(\b _1)_1\gamma ^2_1Y^2\\
&&\hspace*{1cm}\times
\smi (\tqlb \tilde{p}^2_2g^2)_2(e_i)_2P^2)e_j\>\\
&{{\rm (\ref{pql}, \ref{ca})}\atop =}&\sum \limits _{i, j=1}^n
\< e^i, S^{-2}\left(S(\b _2g^2)\gamma ^1Y^1(e_i)_1P^1g^1
S(f^2\gamma ^2_2Y^3)\right)\b _{(1, 1)}\>\\
&&\hspace*{1cm}
\< e^j, \smi (f^1\gamma ^2_1Y^2(e_i)_2P^2)\b _{(1, 2)}e_j\>\\
&{{\rm (\ref{f12})}\atop =}&\sum \limits _{i, j=1}^n
\< e^i, S^{-2}\left(S(Y^3\b _2g^2)f^1\gamma ^1_1(e_i)_1(\b _1)_{(1, 1)}P^1g^1
S(S(Y^1)\gamma ^2)\right)\>
\end{eqnarray*}
\begin{eqnarray*}
&&\hspace*{1cm}
\< e^j, (\b _1)_2\smi (S(Y^2)f^2\gamma ^1_2(e_i)_2(\b _1)_{(1, 2)}P^2)e_j\>\\
&{{\rm (\ref{qr1}, \ref{gdf})}\atop =}&\sum \limits _{i, j=1}^n
\< e^i, \gamma ^1S^{-2}\left(S(Y^3\d ^2)f^1(e_i)_1P^1\d ^1
S(S(Y^1)\gamma ^2)\right)\>\\
&&\hspace*{1cm}
\< e^j, \smi (S(Y^2)f^2(e_i)_2P^2)e_j\>\\
&{{\rm (\ref{gd}, \ref{d2})}\atop =}&\sum \limits _{i, j=1}^n
\< \ov{S}^{-2}(e^i), S(X^1x^1_1Y^1)\a x^3y^3_2Z^3
\smi (f^1P^1y^1\b )(e_i)_2Y^3y^2Z^1\b \\
&&\hspace*{1cm}\times
S(S(X^2)\a X^3x^2y^3_1Z^2)S^2(x^1_2)\>
\< e^j, \smi (f^2P^2)(e_i)_1Y^2e_j\>\\
&{{\rm (\ref{f7}, \ref{qr})}\atop =}&\sum \limits _{i, j=1}^n
\< \ov{S}^{-2}(e^i), S(q^1x^1_1Y^1)\a x^3y^3_2Z^3
\smi (f^1P^1y^1\b)(e_i)_2x^1_{(2, 2)}Y^3\\
&&\hspace*{1cm}\times
y^2Z^1\b S(S(q^2)x^2y^3_1Z^2)\>\< e^j, \smi (f^2P^2)(e_i)_1x^1_{(2, 1)}Y^2e_j\>\\
&{{\rm (\ref{q1}, \ref{f7}, \ref{ca})}\atop =}&\sum \limits _{i, j=1}^n
\< \ov{S}^{-2}(e^i), S(q^1Y^1)\a x^3y^3_2Z^3\smi (f^1(x^1_1)_{(1, 1)}P^1y^1\b)
(e_i)_2Y^3x^1_2\\
&&\hspace*{1cm}\times
y^2Z^1\b S(S(q^2)x^2y^3_1Z^2)\>\\
&&\hspace*{1cm}
\< e^j, (x^1_1)_2\smi (f^2(x^1_1)_{(1, 2)}P^2)(e_i)_1Y^2e_j\>\\
&{{\rm (\ref{qr1}, \ref{q3}, \ref{q5})}\atop =}&\sum \limits _{i, j=1}^n
\< \ov{S}^{-2}(e^i), S(q^1Y^1)\a \smi (f^1P^1p^1\b)(e_i)_2Y^3p^2S^2(q^2)\>\\
&&\hspace*{1cm}
\< e^j, \smi (f^2P^2)(e_i)_1Y^2e_j\>\\
&{{\rm (\ref{f7}, \ref{qr2})}\atop =}&\sum \limits _{i, j=1}^n
\< \ov{S}^{-2}(e^i), S(q^1Q^1_1x^1_{(1, 1)})\a \smi (f^1P^1p^1\b)
(e_i)_2\smi (x^2g^1)\\
&&\hspace*{1cm}\times
Q^2x^1_2p^2\>
\< e^j, \smi (f^2P^2)(e_i)_1\smi (x^3g^2)q^2Q^1_2x^1_{(1, 2)}e_j\>\\
&{{\rm (\ref{ca})}\atop =}&\sum \limits _{i, j=1}^n
\< \ov{S}^{-1}(e^i), S(q^1Q^1_1x^1_{(1, 1)})\a \smi (x^2(e_i)_1P^1p^1\b)Q^2x^1_2p^2\>\\
&&\hspace*{1cm}
\< e^j, \smi (x^3(e_i)_2P^2)q^2Q^1_2x^1_{(1, 2)}e_j\>\\
&{{\rm (\ref{f7})}\atop =}&\sum \limits _{i, j=1}^n
\< \ov{S}^{-1}(e^i), \a \smi \left(x^2(e_i)_1q^1_1(Q^1x^1_1)_{(1, 1)}P^1p^1\b \right)
Q^2x^1_2p^2\>\\
&&\hspace*{1cm}
\< e^j, q^2(Q^1x^1_1)_2\smi \left(x^3(e_i)_2q^1_2(Q^1x^1_1)_{(1, 2)}P^2\right)e_j\>\\
&{{\rm (\ref{qr1}, \ref{pqr})}\atop =}&\sum \limits _{i, j=1}^n
\< \ov{S}^{-1}(e^i), \a \smi (x^2(e_i)_1Q^1x^1_1p^1\b)Q^2x^1_2p^2\>
\< e^j, \smi (x^3(e_i)_2)e_j\>\\
&{{\rm (\ref{f13}, \ref{f7})}\atop =}&\sum \limits _{i, j=1}^n
\< \ov{S}^{-1}(e^i), \a \smi (\tqla (e_i)_1X^1p^1\b )X^2p^2S(X^3)\>
\< e^j, \smi (\tqlb (e_i)_2)e_j\>\\
&{{\rm (\ref{qr}, \ref{f10})}\atop =}&\sum \limits _{j=1}^n
\< e^j, \smi \left(\tilde{\mf{P}}(S^{-2}(\b S(\a)))\right)e_j\>\\
&=&\va \left(\tilde{\mf{P}}(S^{-2}(\b S(\a)))\right)\\
&{{\rm (\ref{f10})}\atop =}&
\sum \limits _{i=1}^n\< e^i, S^{-2}\left(S(\b)\a e_i\b S(\a)\right)\>
={\rm Tr}\left(h\mapsto S^{-2}(S(\b )\a h\b S(\a))\right),
\end{eqnarray*}
where in the last but one equality we used the fact that
$\smi \left(\tilde{\mf{P}}(S^{-2}(\b S(\a)))\right)$ is a right
integral in $H$. So the proof is complete.
\end{proof}
${\;\;\;}$
We will end this Section by computing the representation-theoretic
rank of $D^{\omega}(H)$, the quasi-triangular quasi-Hopf algebra
constructed in \cite{bp}.\\
${\;\;\;}$
Let $H$ be a finite dimensional cocommutative Hopf algebra and
$\omega : H\ot H\ot H\ra k$ a normalized $3$-cocycle on $H$, this means
a convolution invertible map satisfying the conditions:
\begin{eqnarray*}
&&\omega (a_1, b_1, c_1d_1)\omega (a_2b_2, c_2, d_2)=
\omega (b_1, c_1, d_1)\omega (a_1, b_2c_2, d_2)\omega (a_2, b_3, c_3),\\
&&\omega (1, a, b)=\omega (a, 1, b)=\omega (a, b, 1)=\va (a)\va (b),
\end{eqnarray*}
for all $a, b, c, d\in H$. Identifying $(H\ot H\ot H)^*$ with
$H^*\ot H^*\ot H^*$ we can regard $\omega$ and its convolution inverse
$\omega^{-1}$ as elements of $H^*\ot H^*\ot H^*$. Then the commutative Hopf algebra
$H^*$ has a non-trivial quasi-Hopf algebra structure
by keeping the usual multiplication, unit,
comultiplication, counit and antipode of $H^*$, and defining the reassociator
$\Phi =\omega ^{-1}$ and the
elements $\a =\va$, $\b (h)=\omega (h_1, S(h_2), h_3)$, $h\in H$. We shall
denote by $H^*_{\omega}$ the quasi-Hopf algebra structure on $H^*$ defined above.\\
${\;\;\;}$
Now, roughly speaking, the quasi-Hopf algebra $D^{\omega}(H)$ can be identified
as a quasi-triangular quasi-Hopf algbebra
with $D(H^*_{\omega})$, the quantum double associated to the finite
dimensional quasi-Hopf algebra $H^*_{\omega}$. Note that this point of view was
given in \cite{pv}, the initial construction of $D^{\omega}(H)$ being presented
earlier in \cite{bp} as a generalization of the Dijkgraaf-Pasquier-Roche quasi-Hopf 
algebra $D^{\omega}(G)$ constructed in \cite{dpr} (here $G$ is a finite group 
and $\omega$ is a normalized $3$-cocycle on $G$).\\
${\;\;\;}$
Having this description for $D^{\omega}(H)$ and the
result in Proposition \ref{pr4.5} we can easily compute its 
representation-theoretic rank. Note that, one of the goals in \cite{bpv3} was 
to compute this rank but at that moment only a partial answer was given.

\begin{proposition}\label{pr4.6}
Let $H$ be a finite dimensional cocommutative Hopf algebra and
$\omega$ a normalized $3$-cocycle on $H$.
Then $\un{\rm dim}(D^{\omega}(H))={\rm dim}(H)$.
\end{proposition}
\begin{proof}
By Proposition \ref{pr4.5} we have that 
\[
\un{\rm dim}(D^{\omega}(H))={\rm Tr}\left(\v \mapsto
\ov{S}^{-2}(\ov{S}(\b)\a \v \b \ov{S}(\a))\right).
\]
Since $H^*$ is commutative we have that $\ov{S}^2=id_{H^*}$. 
Moreover, $\a=\va$ and from \cite{bp} we know that $\b$ is convolution 
invertible with $\b ^{-1}=\ov{S}(\b )$. 
Therefore, the above formula comes out explicitly as 
\[
\un{\rm dim}(D^{\omega}(H))={\rm Tr}(\v \mapsto \v)={\rm Tr}(id_{H^*})=
{\rm dim}(H^*)={\rm dim}(H),
\]
and this ends the proof.
\end{proof}
\section{The trace formula for quasi-Hopf algebras}\selabel{5}
\setcounter{equation}{0}
${\;\;\;}$
When $H$ is an ordinary Hopf algebra the formula in
Proposition \ref{pr4.5} reduces to
$\un{\rm dim}(H)=\un{\rm dim}(D(H))={\rm Tr}(S^{-2})={\rm Tr}(S^2)$.
As we has already explained in Introduction, using the
Radford and Larson results \cite{lr1, lr2} on one hand, and the Etingof and Gelaki
result \cite{eg} on the other hand, we obtain that
\[
\un{\rm dim}(H)=\un{\rm dim}(D(H))=\left \{\begin{array}{lc}
0&\hspace*{-3mm}
\mbox{,~if $H$ is neither semisimple or cosemisimple}\\
{\rm dim}(H)&\hspace*{-3mm}
\mbox{,~if $H$ is both semisimple and cosemisimple}.
\end{array}\right.
\]
${\;\;\;}$
In this Section we will generalize to the quasi-Hopf algebra setting the first
result. Even if a quasi-Hopf algebra is not a coassociative coalgebra,
as we have seen in Introduction we can define the cosemisimple notion.
Let us explain this more precisely.\\
${\;\;\;}$
Let $H$ be a finite dimensional quasi-Hopf algebra and $t$ a non-zero right
integral in $H$. Since $\int _l$ is a two-sided ideal of $H$,
it follows from the uniqueness of the integrals in $H$ that
there exists $\mu \in H^*$ such that
\[
th=\mu (h)t,~~\forall~~t\in \int _l~~{\rm and}~~h\in H.
\]
Note that $\mu$ is an algebra map;
as in the Hopf case we will call $\mu$ the distinguished group-like
element of $H^*$. We notice that $\mu =\va$ if and only if $H$ is 
unimodular, this means if and only if $\int _l=\int _r$.\\
${\;\;\;}$
Now, following \cite{hn3}, a left cointegral in $H$
is an element $\l \in H^*$ such that
\[
\l (V^2h_2U^2)V^1h_1U^1=\mu (x^1)\l (hS(x^2))x^3,~~
\forall~~h\in H,
\]
where $U=U^1\ot U^2$ is the element defined in (\ref{U}) and
\[
V=V^1\ot V^2:=\smi (f^2p^2)\ot \smi (f^1p^1).
\]
${\;\;\;}$
We will say that a left cointegral $\l$ is normalized if
$\l (\smi (\a)\b)=1$ and we will call a finite dimensional quasi-Hopf algebra
$H$ cosemisimple if $H$ has a normalized left cointegral.\\
${\;\;\;}$
By ${\cal L}$ we denote the space of left cointegral in $H$. Then
the map
\begin{equation}\label{nu}
\nu : {\cal L}\ot H\ra H^*,~~
\nu (\l \ot h)(h')=\l (h'S(h))~~\forall ~~\l \in {\cal L}~~{\rm and}~~
h, h'\in H,
\end{equation}
is an isomorphism of right quasi-Hopf bimodules (the definition of a right
quasi-Hopf $H$-bimodule can be found in \cite{hn3}; roughly speaking it is a right
$H$-comodule within the monoidal category of $H$-bimodules). Here
${\cal L}\ot H$ and $H^*$ are right quasi-Hopf $H$-bimodules via the structures
\begin{equation}\label{f14}
{\cal L}\ot H : \left \{\begin{array}{lc}
h'\cd (\l \ot h)\cd h{''}=\mu(h'_1)\l \ot h'_2hh{''}\\
\l \ot h\mapsto \mu(x^1)\l \ot x^2h_1\ot x^3h_2,
\end{array}\right. 
\end{equation}
\begin{equation}\label{f15}
H^* : \left \{\begin{array}{lc}
\< h'\rightharpoondown \v
\leftharpoondown h{''}, h\>=\< \v , \smi (h')hS(h{''})\>\\
\v \mapsto \sum \limits _{i=1}^ne^i*\v \ot e_i,
\end{array}\right.
\end{equation}
for all $\l \in {\cal L}$, $h, h', h{''}\in H$ and $\v \in H^*$,
where we denoted by $*$ the non-associative multiplication on $H^*$
defined for all $\v , \psi \in H^*$ and $h\in H$ by
\[
\< \v *\psi , h\>:=\< \v , V^1h_1U^1\>\< \psi , V^2h_2U^2\>.
\]
${\;\;\;}$
It follows from the above that ${\rm dim}({\cal L})=1$, and that
for a fixed non-zero left cointegral $\l $ in $H$ the isomorphism
$\nu$ defined in (\ref{nu}) induces a right $H$-linear isomorphism
\[
\tilde{\nu} : H\ra H^*, ~~\tilde{\nu}(h)(h')=\l (h'S(h))) ~~
\forall ~~h, h'\in H.
\]
(Here $H$ and $H^*$ are right $H$-modules via the right regular
representation and $(\v \leftharpoondown h)(h')=\v (h'S(h))$, respectively.)
In particular, there is an unique $r\in H$ such that $\tilde{\nu}(r)=\va$, this means
$\l (hS(r))=\va (h)$, for all $h\in H$. As in the Hopf case we can show that
$r$ is a non-zero integral with the property that $\l (S(r))=1$. Indeed,
the fact that $\tilde{\nu}$ is right $H$-linear implies:
\[
\tilde{\nu}(rh)=\tilde{\nu}(r)\leftharpoondown h=\va \leftharpoondown h=
\va(h)\va =\tilde{\nu }(\va(h)r),
\]
for all $h\in H$. Since $\tilde{\nu}$ is bijective we conclude that
$rh=\va(h)r$, for all $h\in H$, i.e. $r\in \int _r$. Now, $\tilde{\nu}(r)=\va$
implies $\l (hS(r))=\va(h)$ for all $h\in H$, and this is equivalent to
$\l (S(r))=1$.\\
${\;\;\;}$
As we will see the pair $(\l , r)$ described above plays an important role in the
trace formula for quasi-Hopf algebras. In particular,
we will obtain an important result characterizing semisimple cosemisimple
quasi-Hopf algebras in terms of the trace of the "square" of the antipode.
Recall that a semisimple quasi-Hopf algebra is a quasi-Hopf algebra
which is semisimple as an algebra.

\begin{theorem}\label{th5.1}
Let $H$ be a finite dimensional quasi-Hopf algebra, $\mu$ the distinguished
group-like element of $H^*$, $\l$ a non-zero left cointegral
in $H$ and $r$ a right integral in $H$ such that $\l (S(r))=1$. Then:
\begin{itemize}
\item[i)] For any endomorphism $\chi $ of $H$ we have that
\[
{\rm Tr}(\chi )=\mu (q^1_1x^1)\l \left(\chi (q^2x^3r_2p^2)S(q^1_2x^2r_1p^1)\right).
\]
\item[ii)] ${\rm Tr}\left(h\mapsto \b S(\a)S^2(h)S(\b)\a \right)=\va(r)\l (\smi (\a)\b)$.
In particular, $H$ is semi-simple and cosemisimple if and only if
${\rm Tr}\left(h\mapsto \b S(\a)S^2(h)S(\b)\a \right)\not=0$.
\end{itemize}
\end{theorem}
\begin{proof}
For any linear morphism $\chi : H\ra H$ we denote by $\chi ^*: H^*\ra H^*$
the dual morphism of $\chi$. We also denote by $\eta :H^*\ot H\ra {\rm End}(H^*)$
the linear map defined for all $\v , \psi \in H^*$ and $h\in H$ by
\[
\eta (\v \ot h)(\psi )=\psi (h)\v .
\]
Then, exactly as in \cite[Section 7.4]{dnr}, one can easily see that
\begin{eqnarray}
&&\eta (\v \ot h)\circ \chi ^*=\eta (\v \ot \chi (h)),\label{f16}\\
&&{\rm Tr}(\eta (\v \ot h))=\v (h),\label{f17}
\end{eqnarray}
for all $\v \in H^*$, $h\in H$ and $\chi \in {\rm End}(H)$.\\
${\;\;\;}$
i) The fact that $\nu$ is right $H$-colinear shows by
using of (\ref{f14}, \ref{f15}) that
\[
\v (V^1h_1U^1)\l (V^2h_2U^2S(h'))=\mu (x^1)\v (x^3h'_2)\l (hS(x^2h'_1)),
\]
for all $\v \in H^*$ and $h, h'\in H$. If we write the above
equation for $h'=r$ and use the fact that $S(r)\in \int _l$ such that
$\l (S(r))=1$, we obtain
\[
\v (\smi (\b )h\a )=\mu (x^1)\v (x^3r_2)\l (hS(x^2r_1)),
\]
for all $\v \in H^*$ and $h\in H$. In particular, we have that
\begin{eqnarray*}
&&\< p^2\rh \v \lh q^2, \smi (\b)\smi (q^1)hS(p^1)\a \>\\
&&\hspace*{1cm}
=\mu (x^1)\< p^2\rh \v \lh q^2, x^3r_2\>\l (\smi (q^1)hS(p^1)S(x^2r_1)),
\end{eqnarray*}
and this comes out explicitly as
\[
\v (h)=\mu (q^1_1x^1)\v (q^2x^3r_2p^2)\l (hS(q^1_2x^2r_1p^1)),
\]
for all $\v \in H^*$ and $h\in H$, where we used the formula
\begin{equation}\label{f18}
\l (\smi (h)h')=\mu (h_1)\l (h'S(h_2)),~~
\forall ~~h, h'\in H,
\end{equation}
which can be found in \cite[Lemma 3.3]{bc1}. In other words we have obtained
\begin{equation}\label{f19}
\eta (\l \leftharpoondown q^1_2x^2r_1p^1\ot \mu(q^1_1x^1)q^2x^3r_2p^2)=id_{H^*}.
\end{equation}
Now, using (\ref{f16}), (\ref{f17}) and the fact that
${\rm Tr}(\chi )={\rm Tr}(\chi ^*)$ we conclude that
\begin{eqnarray*}
{\rm Tr}(\chi )&=&{\rm Tr}(\chi ^*)={\rm Tr}(id_{H^*}\circ \chi ^*)\\
&=&{\rm Tr}\left(\eta (\l \leftharpoondown q^1_2x^2r_1p^1\ot \mu(q^1_1x^1)q^2x^3r_2p^2)
\circ \chi^*\right)\\
&=&{\rm Tr}\left(\eta (\l \leftharpoondown q^1_2x^2r_1p^1\ot \mu(q^1_1x^1)\chi(q^2x^3r_2p^2)
\right)\\
&=&\mu (q^1_1x^1)\l \left(\chi (q^2x^3r_2p^2)S(q^1_2x^2r_1p^1)\right).
\end{eqnarray*}
${\;\;\;}$
ii) One can easily see that (\ref{pqr}) and $r\in \int _r$ imply:
\[
r_1\ot r_2=r_1q^1_1p^1\ot r_2q^1_2p^2S(q^2)=r_1p^1\ot r_2p^2\a .
\]
Also, by (\ref{qr1}) we have
\[
r_1p^1h\ot r_2p^2=(rh_1)_1p^1\ot (rh_1)_2p^2S(h_2)=r_1p^1\ot r_2p^2S(h),
\]
for any $h\in H$. Combining the two relations above we obtain
\begin{equation}\label{f20}
r_1\ot r_2=r_1p^1\ot r_2p^2\a =r_1p^1\smi (\a )\ot r_2p^2.
\end{equation}
Now, by part i) we have
\begin{eqnarray*}
&&\hspace*{-2.5cm}
{\rm Tr}\left(h\mapsto \b S(\a)S^2(h)S(\b)\a \right)\\
&=&\mu (q^1_1x^1)\l \left(\b S(\a )S^2(q^2x^3r_2p^2)S(\b)\a S(q^1_2x^2r_1p^1)\right)\\
&{{\rm (\ref{f20})}\atop =}&
\mu (q^1_1x^1)\l \left(\b S(\a )S(q^1_2x^2r_1\b S(q^2x^3r_2))\right)\\
&{{\rm (\ref{q5}, \ref{qr})}\atop =}&
\va (r)\mu (q^1_1p^1)\l \left(\b S(\a )S(q^1_2p^2S(q^2))\right)\\
&{{\rm (\ref{pqr})}\atop =}&\va (r)\l (\b S(\a )).
\end{eqnarray*}
Next, we claim that $\va (r)\l (\b S(\a ))=\va (r)\l (\smi (\a )\b )$. Indeed,
if $H$ is not semisimple then by \cite{p} we have that $\va (r)=0$ and
therefore $\va (r)\l (\b S(\a ))=\va (r)\l (\smi (\a )\b )=0$. On the other hand,
if $H$ is semisimple then by the same result in \cite{p} we have that
$\va (\int _l)=\va (\int _r)\not=0$.
In this situation, applying similar arguments as in the Hopf
algebra case we can prove that $H$ is
unimodular, so $\mu =\va$. Finally, by (\ref{f18}) we get
\[
\l (\smi (\a )\b )=\mu (\a _1)\l (\b S(\a _2))=\va (\a _1)\l (\b S(\a _2))
=\l (\b S(\a )),
\]
as claimed. Thus the proof is finished.
\end{proof}
${\;\;\;}$
As a consequence of Proposition \ref{pr4.5} and Theorem \ref{th5.1} we obtain
the following formula for the representation-theoretic ranks of $H$ and $D(H)$.

\begin{theorem}\label{th5.2}
Let $H$ be a finite dimensional quasi-Hopf algebra, $\l $ a left cointegral in
$H$ and $r$ a right integral in $H$ such that $\l (r)=1$. Then
\[
\un{\rm dim}(H)=\un{\rm dim}(D(H))=\va (r)\l (\smi (\a)\b )=
\va _D(\b \rh \l \bowtie r).
\]
In particular, if $H$ is not semisimple or cosemisimple then
\[
\un{\rm dim}(H)=\un{\rm dim}(D(H))=0.
\]
\end{theorem}
\begin{proof}
By $\l _{\rm op}$ we denote a left cointegral in $H^{\rm op}$. It is
straightforward to check that in $H^{\rm op}$ we have
$\mu _{\rm op}=\mu ^{-1}:=\mu \circ S$, and that
the roles of $U$ and $V$ interchange. So $\l _{\rm op}$ is an element of $H^*$
satisfying
\[
\l _{\rm op}(V^2h_2U^2)V^1h_1U^1=\mu ^{-1}(X^1)\l _{\rm op}(\smi (X^2)h)X^3,~~
\forall ~~h\in H.
\]
Note that, if $H$ is unimodular then $\mu =\va$ and therefore a left cointegral
in $H^{\rm op}$ is nothing else than a left cointegral in $H$.\\
${\;\;\;}$
Applying now Theorem \ref{th5.1} to the quasi-Hopf algebra $H^{\rm op}$ we
obtain
\[
{\rm Tr}\left( h\mapsto S^{-2}(S(\b)\a h\b S(\a )\right)=\va (t)\l _{\rm op}(\smi (\a)\b ),
\]
where $t$ is a left integral in $H$ such that $\l _{\rm op}(\smi (t))=1$. If we denote
$r=\smi (t)$ we get that $r$ is a right integral in $H$ such that $\l _{\rm op}(r)=1$.
It follows that $\va (t)=\va (r)$, and that
\[
\un{\rm dim}(H)=\un{\rm dim}(D(H))=
{\rm Tr}\left( h\mapsto S^{-2}(S(\b)\a h\b S(\a )\right)=\va (r)\l _{\rm op}(\smi (\a)\b ).
\]
Finally, we apply the same trick as in the proof of the above Theorem. Namely,
if $H$ is not semisimple then $\va (r)=0$ and we are done. If $H$ is semisimple
then it is unimodular. In this case we have seen that $\l _{\rm op}$ is a cointegral in $H$
and since $\l _{\rm op}(r)=1$ the above equality finishes the proof.
\end{proof}
\begin{remark}
It is conjectured in \cite{hn3} that $\b\rh \l\bowtie r$ is a left integral
in $D(H)$. If it is the case then by the Maschke's theorem proved in \cite{p}
we obtain that $\un{\rm dim}(H)=\un{\rm dim}(D(H))\not=0$ if and only if
$D(H)$ is a semisimple quasi-Hopf algebra. Now, we conjecture that $D(H)$ is
semisimple if and only if $H$ is both semisimple and cosemisimple, if and only
if $h\mapsto S^{-2}(S(\b)\a h\b S(\a ))=id_H$. If it is true then the scalar
$\un{\rm dim}(H)=\un{\rm dim}(D(H))$ has the same value as in the Hopf algebra case.
\end{remark}

\end{document}